\documentclass{amsart}

\usepackage{amssymb}
\usepackage{color}

\allowdisplaybreaks

\numberwithin{equation}{section}

\newtheorem{theorem}{Theorem}[section]
\newtheorem{lemma}[theorem]{Lemma}
\newtheorem{proposition}[theorem]{Proposition}

\newtheorem{remark}[theorem]{Remark}

\newtheorem{TheoA}{Theorem A} \newtheorem{TheoB}{Theorem B}

 \newcommand{\Z}{\mathbb{Z}}
\newcommand{\R}{\mathbb{R}} \newcommand{\C}{\mathbb{C}}
\newcommand{\F}{\mathbb{F}}

\newcommand{\otens}{\overline{\tens}}
\newcommand{\summ}{\sum\nolimits} 
\newcommand{\tens}{\otimes}  \def\G{\mathrm{G}}
\def\1{\mathbf{1}}  
\def\M{\mathcal{M}} \def\A{\mathcal{A}} 
  
 \def\k{\underline{k}} \def\l{\underline{\ell}}

\def\jbar{\underline{j}}
\def\ibar{\underline{i}}

 \newcommand{\dem}{\noindent
  {\bf Proof. }}  

\newcommand{\demA}{\noindent {\bf Proof of Theorem A. }}
  
\newcommand{\demB}{\noindent {\bf Proof of Theorem B. }}

\newcommand{\fin}{\hspace*{\fill} $\square$ \vskip0.2cm}

\newcommand{\si}{\sigma} \newcommand{\eps}{\varepsilon}

\newcommand{\Acirc}{\stackrel{\circ}{A_j}}
\newcommand{\Hcirc}{\stackrel{\circ}{\mathcal{H}_j}}

\addtocounter{tocdepth}{-1}

\begin{document}

\addtolength{\parskip}{+1ex}

\title[Hypercontractivity for free products]{Hypercontractivity for
free products}

\author[Junge, Palazuelos, Parcet, Perrin and Ricard] {Marius Junge, Carlos
  Palazuelos, \\ Javier Parcet, Mathilde Perrin and \'Eric Ricard}

\maketitle

\begin{abstract}
In this paper, we obtain optimal time hypercontractivity bounds for
the free product extension of the Ornstein-Uhlenbeck semigroup acting
on the Clifford algebra. Our approach is based on a central limit
theorem for free products of spin matrix algebras with mixed
commutation/anticommutation relations. With another use of Speicher's
central limit theorem, we may also obtain the same bounds for free
products of $q$-deformed von Neumann algebras interpolating between
the fermonic and bosonic frameworks. This generalizes the work of
Nelson, Gross, Carlen/Lieb and Biane. Our main application yields
hypercontractivity bounds for the free Poisson semigroup acting on the
group algebra of the free group $\F_n$, uniformly in the number of
generators.
\end{abstract}

\section*{{\bf Introduction}}

The two-point inequality was first proved by Bonami and rediscovered
years later by Gross \cite{Bonami,Gross1}. In the context of
harmonic analysis, this inequality was central for Bonami's work on
the relation between integrability of a function and the decay
properties of its Fourier coefficients. It was also instrumental in
Beckner's theorem on the optimal constants for the Hausdorff-Young
inequality \cite{Beckner}. On the other hand, motivated by quantum field theory,
Gross used it as a key step towards his logarithmic Sobolev
inequalities \cite{Gross1}. More recently, the two-point inequality
has also produced very important applications in computer science and
in both classical and quantum information theory
\cite{BRSW,GKKRW,KhVi,KlRe}. If $1 < p \le q < \infty$ and
$\alpha,\beta \in \C$, Bonami-Gross inequality can be
rephrased for $r=e^{-t}$ as follows $$\Big( \sum_{\varepsilon = \pm 1}
\Big| \frac{(1+ \varepsilon r) \alpha + (1-\varepsilon r)
  \beta}{2^{1+\frac{1}{q}}} \Big|^q \Big)^{\frac1q} \le \Big(
\frac{|\alpha|^p + |\beta|^p}{2} \Big)^{\frac1p} \ \Leftrightarrow \ r
\le \sqrt{\frac{p-1}{q-1}}.$$ It can be regarded ---from Bonami's
viewpoint--- as the optimal hypercontractivity bound for the \lq\lq
Poisson semigroup" on the group $\Z_2$, while Gross understood it as
the optimal hypercontractivity bound for the Ornstein-Uhlenbeck
semigroup on the Clifford algebra with one generator
$\mathcal{C}(\R)$. Although the two-point inequality can be
generalized in both directions, harmonic analysis has developed
towards other related norm inequalities in the classical groups
---like $\Lambda_p$ sets in $\Z$--- instead of analyzing the
hypercontractivity phenomenon over the compact dual of other discrete
groups. Namely, to the best of our knowledge only the cartesian
products of $\Z_2$ and $\Z$ have been understood so far, see
\cite{Weissler}. The first goal of this paper is to replace cartesian
products by free products, and thereby obtain hypercontractivity
inequalities for the free Poisson semigroups acting on the group von
Neumann algebras associated to $\F_n = \Z * \Z * \cdots * \Z$ and
$\mathbb{G}_n = \Z_2 * \Z_2 * \cdots * \Z_2$.

Let $\G$ denote any of the free products considered above and let
$\lambda: \G \to \mathcal{B}(\ell_2(\G))$ stand for the corresponding
left regular representation. The group von Neumann algebra
$\mathcal{L}(\G)$ is the weak operator closure of the linear span of
$\lambda(\G)$. If $e$ denotes the identity element of $\G$, the
algebra $\mathcal{L}(\G)$ comes equipped with the standard trace
$\tau(f) = \langle \delta_e, f \delta_e \rangle$. Let
$L_p(\mathcal{L}(\mathrm{G}), \tau)$ be the $L_p$ space over the
noncommutative measure space $(\mathcal{L}(\mathrm{G}), \tau_\G)$
---the so called noncommutative $L_p$ spaces--- with norm $\|f\|_p^p =
\tau |f|^p$. We invite the reader to check that
$L_p(\mathcal{L}(\G),\tau) = L_p(\mathbb{T})$ for $\G = \Z$ after
identifying $\lambda_{\Z}(k)$ with $e^{2\pi i k \cdot}$. In the
general case, the absolute value and the power $p$ are obtained from
functional calculus for this (unbounded) operator on the Hilbert space
$\ell_2(\G)$, see \cite{PX} for details. If $f = \sum_g \widehat{f}(g)
\lambda(g)$, the free Poisson semigroup on $\G$ is given by the family
of linear maps $$\mathcal{P}_{\G,t} f \, = \, \sum_{g \in \G} e^{-t
  |g|} \widehat{f}(g) \lambda(g) \quad \mbox{with} \quad t \in \R_+.$$
In both cases $\G \in \{ \F_n,\mathbb{G}_n \}$, $|g|$ refers to the
Cayley graph length. In other words, $|g|$ is the number of letters
(generators and their inverses) which appear in $g$ when it is written
in reduced form. It is known from \cite{Ha} that $\mathcal{P}_\G =
(\mathcal{P}_{\G,t})_{t \ge 0}$ defines a Markovian semigroup of
self-adjoint, completely positive, unital maps on
$\mathcal{L}(\G)$. In particular, $\mathcal{P}_{\G,t}$ defines a
contraction on $L_p(\mathcal L(\G))$ for every $1\leq p\leq
\infty$. The hypercontractivity problem for $1< p \le q < \infty$
consists in determining the optimal time $t_{p,q} > 0$ above
which $$\|\mathcal{P}_{\G,t} f \|_q \, \le \, \|f\|_p \qquad \mbox{for
  all} \qquad t \ge t_{p,q}.$$ In our first result we provide new
hypercontractivity bounds for the free Poisson semigroups on those
group von Neumann algebras. If $g_1, g_2, \ldots, g_n$ stand for the
free generators of $\F_n$, we will also consider the symmetric
subalgebra $\mathcal{A}_{sym}^n$ of $\mathcal{L}(\F_n)$ generated by
the self-adjoint operators $\lambda(g_j) + \lambda(g_j)^*$. In other
words, we set $$\mathcal{A}_{sym}^n \, = \, \big\langle \lambda(g_1) +
\lambda(g_1)^*, \ldots, \lambda(g_n) + \lambda(g_n)^* \big\rangle''.$$

\begin{TheoA}
If $1 < p \le q < \infty$, we find\hskip1pt$:$
\begin{enumerate}
\item[i)] Optimal time hypercontractivity for $\mathbb{G}_n$ $$\big\|
  \mathcal{P}_{\mathbb{G}_n,t}: L_p(\mathcal{L}(\mathbb{G}_n)) \to
  L_q(\mathcal{L}(\mathbb{G}_n)) \big\| \, = \, 1 \ \Leftrightarrow
  \ t \ge \frac12 \log \frac{q-1}{p-1}.$$

\item[ii)] Hypercontractivity for $\F_n$ over twice the optimal
  time $$\big\| \mathcal{P}_{\F_n,t} \hskip1pt : L_p( \hskip0.5pt
  \mathcal{L}(\F_n) \hskip0.5pt ) \to L_q( \hskip0.5pt
  \mathcal{L}(\F_n) \hskip0.5pt ) \hskip1pt \big\| \, = \, 1
  \hskip13pt \mbox{if} \hskip13pt t \ge \log \frac{q-1}{p-1}.$$

\item[iii)] Optimal time hypercontractivity in the symmetric algebra
  $\mathcal{A}_{sym}^n$ $$\big\| \hskip1pt \mathcal{P}_{\F_n,t}
  \hskip1pt : L_p( \hskip1pt \mathcal{A}_{sym}^n \hskip1pt ) \hskip1pt
  \to \hskip1pt L_q( \hskip1pt \mathcal{A}_{sym}^n \hskip1pt )
  \hskip1pt \big\| \, = \, 1 \ \Leftrightarrow \ t \ge \frac12 \log
  \frac{q-1}{p-1}.$$
\end{enumerate}
\end{TheoA}

Theorem A i) extends Bonami's theorem for $\Z_2^n$ to the free product
case with optimal time estimates. According to the applications in
complexity theory and quantum information of Bonami's result, it is
conceivable that Theorem A could be of independent interest in those
areas. These potential applications will be explored in further
research. Theorem A ii) gives the first hypercontractivity estimate
for the free Poisson semigroup on $\F_n$, where a factor 2 is lost
from the expected optimal time. This is related to our probabilistic
approach to the problem and a little distortion must be done to make
$\F_n$ fit in. Theorem A iii) refines this, providing optimal time
estimates in the symmetric algebra $\mathcal{A}_{sym}^n$. We also
obtain optimal time $L_p \to L_2$ hypercontractive estimates for
linear combinations of words with length less than or equal to
one. Apparently, our probabilistic approach in this paper is limited
to go beyond the constant $2$ in the general case. We managed to push it to
$1+\frac 14 \log2 \sim 1.173$ in the last section. Actually, we have recently
found in \cite{JPPP} an alternative combinatorial/numerical
method which yields optimal $L_2 \to L_q$ estimates for $q \in 2\Z$ and
also reduces the general constant to $\log \, 3 \sim 1.099$ for $1 < p \le q <
\infty$. The drawback of this method is the numerical part: the larger
is the number of generators $n$, the harder is to implement and test
certain pathological terms in a computer. In this respect, Theorem A
ii) is complementary since ---at the price of a worse constant--- we
obtain uniform estimates in $n$.

As we have already mentioned, it is interesting to understand the
two-point inequality as the convergence between the
\emph{trigonometric point of view} outlined above and the
\emph{gaussian point of view}, which was developed along the extensive
study of hypercontractivity carried out in the context of quantum
mechanics and operator algebras. The study of hypercontractivity in
quantum mechanics dates back to the work of Nelson \cite{Nelson1} who
showed that semiboundedness of certain Hamiltonians $H$ associated to
a bosonic system can be obtained from the (hyper)contractivity of the
semigroup $e^{-tA_\gamma}:L_2(\R^d,\gamma)\rightarrow
L_2(\R^d,\gamma)$, where $A_\gamma$ is the Dirichlet form operator for
the Gaussian measure $\gamma$ on $\R^d$. After some contributions
\cite{Glimm,HoSi,Segal} Nelson finally proved in \cite{Nelson2} that
the previous semigroup is contractive from $L_p(\R^d,\gamma)$ to
$L_q(\R^d,\gamma)$ if and only if $e^{-2t}\leq\frac{p-1}{q-1}$; thus
obtaining the same optimal time as in the two-point inequality. By
that time a new deep connection was shown by Gross in \cite{Gross1},
who established the equivalence between the hypercontractivity of the
semigroup $e^{-tA_\mu}$, where $A_\mu$ is the Dirichlet form operator
associated to the measure $\mu$, and the logarithmic Sobolev
inequality verified by $\mu$. During the next 30 years
hypercontractivity and its equivalent formulation in terms of
logarithmic Sobolev inequalities have found applications in many
different areas of mathematics like probability theory, statistical
mechanics or differential geometry. We refer the survey \cite{Gross3}
for an excellent exposition of the topic.

The extension of Nelson's theorem to the fermonic case started with
Gross' papers \cite{Gross4,Gross2}.  Namely, he adapted the argument in
the bosonic case by considering a suitable Clifford algebra $\mathcal
C(\R^d)$ on the fermion Fock space and noncommutative $L_p$ spaces on
this algebra after Segal \cite{Segal2}. In particular,
hypercontractivity makes perfectly sense in this context by
considering the corresponding Ornstein-Uhlenbeck
semigroup $$\mathcal{O}_t \, := \, e^{-tN_0}: \, L_2(\mathcal
C(\R^d),\tau) \rightarrow L_2(\mathcal C(\R^d),\tau).$$ Here $N_0$
denotes the fermion number operator, see Section \ref{Section:
  preliminars} for the construction of the Clifford algebra $\mathcal
C(\R^d)$ and a precise definition of the Ornstein-Uhlenbeck semigroup
on fermion algebras. After some partial results
\cite{Gross2,Lindsay,LiMe}, the optimal time hypercontractivity bound
in the fermionic case was finally obtained by Carlen and Lieb in
\cite{CaLi} $$\big\| \mathcal{O}_t: L_p(\mathcal C(\R^d)) \rightarrow
L_q(\mathcal C(\R^d)) \big\| \, = \, 1 \ \Leftrightarrow \ t \ge
\frac12 \log \frac{q-1}{p-1}.$$ The proof deeply relies on the optimal
2-uniform convexity for matrices from \cite{BCL}.

Beyond its own interest in quantum mechanics, these contributions
represent the starting point of hypercontractivity in the
noncommutative context. This line was continued by Biane \cite{Biane},
who extended Carlen and Lieb's work and obtained optimal time
estimates for the $q$-Gaussian von Neumann algebras $\Gamma_q$
introduced by Bozejko, K\"ummerer and Speicher \cite{BKS}. These
algebras interpolate between the bosonic and fermonic frameworks,
corresponding to $q = \pm 1$. The semigroup for $q=0$ acts diagonally
on free semi-circular variables ---instead of free generators as in
the case of the free Poisson semigroup--- in the context of
Voiculescu's free probability theory \cite{VDN}. We also refer to
\cite{Janson,Kemp,K1,K2,LeRi} for related results in this line. On the
other hand, the usefulness of the two-point inequality in the context
of computer science has motivated some other extensions to the
noncommutative setting more focused on its applications to quantum
computation and quantum information theory. In \cite{BRW}, the authors
studied extensions of Bonami's result to matrix-valued functions $f:
\Z_2^n \rightarrow M_n(\C)$, finding optimal estimates for $q=2$ and
showing some applications to coding theory. In \cite{MoOs}, the
authors introduced quantum boolean functions and obtained
hypercontractivity estimates in this context with some consequences in
quantum information theory, see also the recent work \cite{Montanaro}.

The very nice point here is that, although our main motivation to
study the Poisson semigroup comes from harmonic analysis, we realized
that a natural way to tackle this problem is by means of studying the
Ornstein-Uhlenbeck semigroup on certain von Neumann algebras. In
particular, a significant portion of Theorem A follows from our main
result, which extends Carlen and Lieb's theorem to the case of free
product of Clifford algebras. The precise definitions of reduced free
products which appear in the statement will be recalled for the
non-expert reader in the body of the paper.
\begin{TheoB}
Let $\M_\alpha = \mathcal{C}(\R^{d_\alpha})$ be the Clifford algebra
with $d_\alpha$ generators for any $1 \le \alpha \le n$ and construct
the corresponding reduced free product von Neumann algebra $\M = \M_1
* \M_2 * \cdots * \M_n$. If $\mathcal{O}_\alpha =
(\mathcal{O}_{\alpha,t})_{t \ge 0}$ denotes the Ornstein-Uhlenbeck
semigroup acting on $\M_\alpha$, consider the free product semigroup
$\mathcal{O}_\mathcal{M} = (\mathcal{O}_{\M,t})_{t \ge 0}$ given by
$\mathcal{O}_{\M,t} = \mathcal{O}_{1,t} * \mathcal{O}_{2,t} * \cdots *
\mathcal{O}_{n,t}$.  Then, we find for $1 < p \le q < \infty$
$$\big\| \mathcal{O}_{\M,t}: L_p(\M) \to L_q(\M) \big\| \, = \, 1
\quad \Leftrightarrow \quad t \, \ge \, \frac12 \log
\frac{q-1}{p-1}.$$
\end{TheoB}
It is relevant to point out a crucial difference between our approach
and the one followed in \cite{Bonami,CaLi,Nelson2}. Indeed, in all
those cases the key point in the argument is certain basic inequality
---like Bonami's two-point inequality or Ball/Carlen/Lieb's convexity
inequality for matrices--- and the general result follows from an
inductive argument due to the tensor product structure of the
problem. However, no tensor product structure can be found in our
setting (Theorems A and B). In order to face this problem, Biane
showed in \cite{Biane} that certain optimal hypercontractive estimates
hold in the case of spin matrix algebras with mixed
commutation/anticommutation relations, and then applied Speicher's
central limit theorem \cite{Speicher}. In this paper we will extend
Biane's and Speicher's results by showing that a wide range of von
Neumann algebras can also be approximated by these spin
systems. Namely, the proof of Theorem B will show that the same result
can be stated in a much more general context. As we shall explain, we
may consider the free product of Biane's mixed spins algebras which in
turn gives optimal hypercontractivity estimates for the free products
of $q$-deformed algebras with $q_1, q_2, \ldots, q_n \in [-1,1]$.

\vskip5pt


\noindent \textbf{Acknowledgements.} The authors were supported by
ICMAT Severo Ochoa Grant SEV-2011-0087 (Spain); ERC Grant
StG-256997-CZOSQP (EU); NSF Grant DMS-0901457 (USA); MINECO Grants
MTM-2010-16518 \& MTM-2011-26912 and ``Juan de la Cierva'' program
(Spain); and ANR-2011-BS01-008-11 (France).

\section{{\bf Preliminaries}} \label{Section: preliminars}

In this section we briefly review the definition of the CAR algebra
and the Ornstein-Uhlenbeck semigroup acting on it. We also recall the
construction of the reduced free product of a family of von Neumann
algebras and introduce the free Ornstein-Uhlenbeck semigroup on a
reduced free product of Clifford algebras.

\subsection{The Ornstein-Uhlenbeck semigroup}

The standard way to construct a system of $d$ fermion degrees of
freedom is by means of the antisymmetric Fock space. Let us consider
the $d$-dimensional real Hilbert space $\mathcal H_\R=\R^d$ and its
complexification $\mathcal H_\C=\C^d$. Define the Fock space
$$\mathcal F(\mathcal H_\R) \, = \, \C \Omega \oplus
\bigoplus_{m=1}^\infty \mathcal H_\C^{\otimes_m}$$ for some fixed unit
vector $\Omega \in \mathcal H_\C$ called the vacuum. If $S_m$ denotes
the symmetric group of permutations over $\{1,2,\ldots,m\}$ and
$i(\beta)$ the number of inversions of the permutation $\beta$, we
define the hermitian form $\langle \cdot, \cdot \rangle$ on $\mathcal
F(\mathcal H_\R)$ by $\langle \Omega,\Omega \rangle =1$ and the
following identity $$\big\langle f_1\otimes\cdots\otimes f_m ,
g_1\otimes\cdots\otimes g_n \big\rangle \, = \, \delta_{mn}
\sum_{\beta \in S_m}^{\null} (-1)^{i(\beta)} \langle
f_1,g_{\beta(1)}\rangle \cdots \langle f_m, g_{\beta(m)} \rangle.$$ It
is not difficult to see that the hermitian form $\langle \cdot, \cdot
\rangle$ is non-negative. Therefore, if we consider the completion of
the quotient by the corresponding kernel, we obtain a Hilbert space
that we will call again $\mathcal F(\mathcal H_\R)$. Let us denote by
$(e_j)_{j=1}^d$ the canonical basis of $\mathcal H_\R = \R^d$. Then,
we define the $j$-th fermion annihilation operator acting on $\mathcal
F(\mathcal H_\R)$ by linearity as $c_j(\Omega)=0$
and $$c_j(f_1\otimes\cdots \otimes f_m) \, = \, \sum_{i=1}^m(-1)^{i-1}
\langle f_i,e_j\rangle \, f_1 \otimes\cdots\otimes f_{i-1}\otimes
f_{i+1}\otimes\cdots\otimes f_m.$$ Its adjoint $c_j^*$ is called the
$j$-th fermion creation operator on $\mathcal F(\mathcal H_\R)$. It is
determined by $c_j^*(\Omega)=e_j$ and $c_j^*(f_1\otimes\cdots \otimes
f_m)=e_j\otimes f_1\otimes\cdots\otimes f_m$. It is quite instrumental
to observe that $c_ic_j + c_jc_i = 0$ and $c_i c_j^* + c_j^*c_i =
\delta_{ij} \mathbf{1}$. The basic free Hamiltonian on $\mathcal
F(\mathcal H_\R)$ is the fermion number operator $$N_0 = \sum_{j=1}^d
c_j^*c_j.$$ It generates the fermion oscillator semigroup $(
\exp(-tN_0) )_{t \ge 0}$. Then, one defines the configuration
operators $x_j=c_j+c_j^*$ for $1 \le j \le d$. Denote by $\mathcal
C(\R^d)$ the unit algebra generated by them. Note that these operators
verify the canonical anti-commutation relations
(CAR) $$x_ix_j+x_jx_i=2\delta_{ij} \qquad \mbox{ and } \qquad
x_j^*=x_j.$$ It is well-known that $\mathcal C(\R^d)$ can be
concretely represented as a subalgebra of the matrix algebra $\mathbb{M}_{2^d}$
by considering $d$-chains formed by tensor products of Pauli
matrices. The key point for us is that the $2^d$ distinct monomials in
the $x_j$'s define a basis of $\mathcal C(\R^d)$ as a vector
space. Indeed, given any subset $A$ of $[d] := \{1,2,\ldots,d\}$ we
shall write $x_A = x_{j_1} x_{j_2} \cdots x_{j_s}$ where $(j_1, j_2,
\ldots, j_s)$ is an enumeration of $A$ in increasing order. If we also
set $x_\emptyset = \mathbf{1}$, it turns out that $\{ x_A \, | \ A
\subset [d] \}$ is a linear basis of $\mathcal C(\R^d)$. In
particular, any $X \in \mathcal C(\R^d)$ has the form
$$X \, = \, \alpha_\emptyset \mathbf{1} \, + \, \sum_{s=1}^d \, \sum_{1 \le j_1
  < \cdots < j_s \le d} \alpha_{j_1, \ldots, j_s} x_{j_1} \cdots
x_{j_s}.$$
The vacuum $\Omega$ defines a tracial state $\tau$ on $\mathcal
C(\R^d)$ by $\tau(X)=\langle X\Omega,\Omega\rangle$. We denote by
$L_p(\mathcal C(\R^d),\tau)$ or just $L_p(\mathcal C(\R^d))$ the
associated non-commutative $L_p$-space. The map $X\mapsto X\Omega$
defines a continuous embedding of $\mathcal C(\R^d)$ into $\mathcal
F(\R^d)$ which extends to a unitary isomorphism $L_2(\mathcal C(\R^d))
\simeq \mathcal F(\R^d)$. Then, instead of working on the Fock space
$\mathcal F(\R^d)$ and with the semigroup $\exp(-tN_0)$, we can
equivalently consider $\mathcal C(\R^d)$ and the Ornstein-Uhlenbeck
semigroup on $\mathcal C(\R^d)$ defined by
\begin{align*}
\mathcal{O}_t(X) \, = \, \alpha_\emptyset \mathbf{1} \, + \, \sum_{s=1}^d
e^{-ts} \, \sum_{1 \le j_1 < \cdots < j_s \le d} \alpha_{j_1, \ldots,
  j_s} x_{j_1} \cdots x_{j_s}.
\end{align*}
If $1 < p \le q< \infty$, the main result in \cite{CaLi}
yields $$\big\| \mathcal{O}_t: L_p(C(\R^d))\rightarrow L_q(C(\R^d))
\big\| \, = \, 1 \ \Leftrightarrow \ t \ge \frac12 \log
\frac{q-1}{p-1}.$$

\subsection{Free product of von Neumann algebras}

Let $(A_j, \phi_j)_{j \in J}$ be a family of unital C$^*$-algebras
with distinguished states $\phi_j$ whose GNS constructions $(\pi_j,
\mathcal H_j, \xi_j)$ with $\mathcal H_j=L_2(A_j, \phi_j)$ are
faithful. Let us define $$\Acirc \, = \, \Big\{ a \in A_j \, \big|
\ \phi_j(a)=0 \Big\} \qquad \mbox{and} \qquad \Hcirc=\xi_j^\perp$$ so
that $A_j=\mathbb{C} \mathbf{1} \oplus \Acirc$ and $\mathcal
H_j=\mathbb{C}\xi_j \oplus \Hcirc$. Note that we have natural maps
$i_j=A_j \rightarrow \mathcal H_j$ such that $\phi_j(a^*b)=\langle
i_j(a),i_j(b) \rangle_{\mathcal H_j}$ for every $j \in J$. Let us
consider the full Fock space associated to the free
product \vskip-6pt $$\mathcal F \ = \ \mathbb{C}\Omega \, \oplus
\bigoplus_{\substack{m \ge 1\\ j_1 \neq j_2 \neq \cdots \neq j_m}}
\stackrel{\circ}{\mathcal{H}}_{j_1} \otimes\cdots\otimes
\stackrel{\circ}{\mathcal{H}}_{j_m}$$ with inner product $$\big\langle
h_1 \otimes \cdots \otimes h_m, h_1' \otimes \cdots \otimes h_n'
\big\rangle \, = \, \delta_{mn} \prod_{j=1}^m \langle h_j,h_j'
\rangle.$$ Each algebra $A_j$ acts non-degenerately on $\mathcal F$
via the map $\omega_j: A_j\rightarrow \mathcal{B}(\mathcal F)$ in the
following manner. Since we can decompose every $z \in A_j$ as $z =
\phi_j(z) \mathbf{1}+ a$ with $\phi_j(a)=0$, it suffices to define
$\omega_j(a)$. Let $h_1 \otimes \cdots \otimes h_m$ be a generic
element in $\mathcal{F}$ with $h_i \in \mathcal{H}_{j_i} \ominus \C
\xi_{j_i}$. If $j \neq j_1$, we set $$\omega_j(a) \big( h_1\otimes
\cdots \otimes h_m \big) \, = \, i_j(a) \otimes h_1\otimes \cdots
\otimes h_m.$$ When $j = j_1$ we add and subtract the mean to obtain
\begin{eqnarray*}
\omega_j(a) \big( h_1\otimes\cdots \otimes h_m \big) & = & \big\langle
\xi_j, \pi_j(a)(h_1) \big\rangle_{\mathcal{H}_j}
h_2\otimes\cdots\otimes h_m \\ & + & \Big( \pi_j(a)(h_1) - \big\langle
\xi_j,\pi_j(a)(h_1) \big\rangle_{\mathcal{H}_j} \xi_j \Big) \otimes
h_2\otimes\cdots \otimes h_m.
\end{eqnarray*}
The faithfulness of the GNS construction of $(A_j, \phi_j)$ implies
that the representation $\omega_j$ is faithful for every $j \in
J$. Thus, we may find a copy of the algebraic free product $$A \ =
\ \mathbb{C}\Omega \, \oplus \bigoplus_{\substack{m \geq 1\\ j_1\neq
    j_2\neq\cdots\neq j_m}}
\stackrel{\circ}{A_{j_1}}\otimes\cdots\otimes
\stackrel{\circ}{A_{j_m}}$$ in $\mathcal{B}(\mathcal F)$. The reduced
free product of the family $(A_j, \phi_j)_{j \in J}$ is the
C$^*$-algebra generated by these actions. In other words, the norm
closure of $A$ in $\mathcal{B}(\mathcal F)$. It is denoted
by $$(A,\phi) \, = \, *_{j \in J} (A_j, \phi_j),$$ where the state
$\phi$ on $A$ is given by $$\phi(\mathbf{1}) \, = \, 1 \qquad
\mbox{and} \qquad \phi(a_1\otimes\cdots\otimes a_m)=0$$ for $m\geq 1$
and $a_i \in \stackrel{\circ}{A_{j_i}}$ with $j_1\neq
j_2\neq\cdots\neq j_m$. Each $A_j$ is naturally considered as a
subalgebra of $A$ and the restriction of $\phi$ to $A_j$ coincides
with $\phi_j$. It is helpful to think of the elementary tensors above
$a_1\otimes\cdots\otimes a_m$ as words of length $m$, where the empty
word $\Omega$ has length $0$. In this sense, a word
$a_1\otimes\cdots\otimes a_m$ can be identified with the product
$a_1a_2\cdots a_m$ via the formula $a_1\cdots a_m\Omega=a_1\otimes
\cdots\otimes a_m$.

This construction also holds in the category of von Neumann
algebras. Let $(\mathcal M_j, \phi_j)_{j \in J}$ be a family of von
Neumann algebras with distinguished states $\phi_j$ whose GNS
constructions $(\pi_j, \mathcal{H}_j, \xi_j)$ are faithful. Then, the
corresponding reduced free product von Neumann algebra is the weak-$*$
closure of $*_{j \in J}(\mathcal M_j, \phi_j)$ in
$\mathcal{B}(\mathcal F)$ which will be denoted by $(\M,\phi) \, = \,
\overline{*}_{j \in J}(\mathcal M_j,\phi_j).$ As before, the $\mathcal
M_j$'s are regarded as von Neumann subalgebras of $\M$ and the
restriction of $\phi$ to $\mathcal M_j$ coincides with $\phi_j$. A
more complete explanation of the reduced free product of von Neumann
algebras can be found in \cite{VDN}. Let us now consider a family
$(\Lambda_j: \mathcal M_j \rightarrow \mathcal M_j)_{j \in J}$ of
normal, completely positive, unital and trace preserving maps. Then,
it is known from \cite[Theorem 3.8]{BlDy} that there exists a map
$\Lambda \, = \, *_{j \in J} \Lambda_j:\M\rightarrow \M$ such that
$\Lambda (x_1 x_2 \cdots x_m) = \Lambda_{j_1}(x_1) \cdots
\Lambda_{j_m}(x_m)$, whenever $x_i \in \M_{j_i}$ is trace 0 and $j_i
\neq j_{i+1}$ for $1 \le i \le m-1$. This map is called the free
product map of the $\Lambda_j$'s. In particular we may take $\M_j =
\mathcal C(\R^d)$ for $1 \le j \le n$ and $\Lambda_j =
\mathcal{O}_{j,t}$, the Ornstein-Uhlenbeck semigroup on $\M_j$ at time
$t$. The resulting free product maps $\mathcal{O}_{\M} =
(\mathcal{O}_{\M,t})_{t \ge 0}$ with $\mathcal{O}_{\M,t} =
\mathcal{O}_{1,t} * \mathcal{O}_{2,t} * \cdots * \mathcal{O}_{n,t}$
will be referred to as the \emph{free Ornstein-Uhlenbeck semigroup} on
the reduced free product von Neumann algebra $\M$.

\section{{\bf The free Ornstein-Uhlenbeck semigroup}}\label{Section: proof of theorem B}

This section is devoted to the proof of Theorem B. Of course, we may and will
assume for simplicity that $d_\alpha = d$ for all $1 \le \alpha \le
n$. The key idea is to describe the free product of fermion algebras
and the corresponding Ornstein-Uhlenbeck semigroup as the limit
objects of certain spin matrix models and certain semigroups defined
on them. In this sense, we will extend Biane's results \cite{Biane} by
showing that these matrix models can be used to describe a wide range
of operator algebra frameworks.

Note that the free Ornstein-Uhlenbeck semigroup restricted to a single
free copy $\M_\alpha$ coincides with the fermion oscillator semigroup
on $\M_\alpha$. In particular, we know from Carlen and Lieb's theorem
\cite{CaLi} that the optimal time in Theorem B must be greater than or
equal to $\frac12 \log \frac{q-1}{p-1}$. This proves the necessity, it
remains to prove the sufficiency. Given $1 \le \alpha \le n$ and
recalling that $[d]$ stands for $\{1,2,\ldots,d\}$, we denote by
$(x_i^\alpha)_{i \in [d]}$ the generators of $\M_\alpha=\mathcal
C(\R^d)$. A reduced word in the free product $\M = \M_1 * \M_2 *
\cdots * \M_n$ is then of the form
\begin{equation} \label{general element free productI}
x \, = \, x_{A_1}^{\alpha_1}\cdots x_{A_\ell}^{\alpha_\ell}
\end{equation}
with $A_j \subset [d]$ and $\alpha_j \neq \alpha_{j+1}$. The case
$\ell = 0$ refers to the empty word $\mathbf{1}$. If we set $s_j =
|A_j|$ and write $A_j = \{i_{s_1+ \cdots +
  s_{j-1}+1},\ldots,i_{s_1+\cdots+s_{j-1}+s_j}\}$ ---labeling the
indices in a strictly increasing order--- $x$ can be written as
follows
\begin{align}\label{general element free productII}
x \ = \ \overbrace{x_{i_1}^{\alpha_1}\cdots
  x_{i_{s_1}}^{\alpha_1}}^{x_{A_1}^{\alpha_1}} \,
\overbrace{x_{i_{s_1+1}}^{\alpha_2}\cdots
  x_{i_{s_1+s_2}}^{\alpha_2}}^{x_{A_2}^{\alpha_2}} \, \cdots \,
\overbrace{x_{i_{s_1+\cdots +s_{\ell-1}+1}}^{\alpha_\ell}\cdots
  x_{i_{s_1+\cdots
      +s_\ell}}^{\alpha_\ell}}^{x_{A_\ell}^{\alpha_\ell}}.
\end{align}
In what follows, we will use the notation
$|x|=|A_1|+\cdots+|A_\ell|=s_1+\cdots+s_\ell$.

\subsection{Spin matrix model}\label{Subsection: matrix model}

Given $m \geq 1$, we will describe a spin system with mixed
commutation and anticommutation relations which approximates the free
product of fermions $\M$ as $m \to \infty$. Let us first recall the
construction of a spin algebra in general. In our setting, we will
need to consider three indices. This is why we introduce the sets
$\Upsilon=[n]\times [d]\times \mathbb{Z}_+$ and $\Upsilon_m =
[n]\times [d]\times [m]$ for $m\geq 1$.  Let $\varepsilon:
\Upsilon\times \Upsilon\rightarrow \{-1,1\}$ be any map satisfying
\begin{itemize}
\item $\varepsilon$ is symmetric: $\varepsilon(x,y) = \varepsilon
  (y,x)$,

\item $\varepsilon \equiv -1$ on the diagonal: $\varepsilon (x,x) =
  -1$.
\end{itemize}
Given $m \geq 1$, we will write $\varepsilon_m$ to denote the
truncation of $\varepsilon$ to $\Upsilon_m \times
\Upsilon_m$. Consider the complex unital algebra $\mathcal
A_{\varepsilon_m}$ generated by the elements
$(x_i^\alpha(k))_{(\alpha,i,k)\in \Upsilon_m}$ which satisfy the
commutation/anticommutation relations
\begin{align}\label{relations free}
x_i^\alpha(k)x_j^\beta(\ell)-\varepsilon\big((\alpha,i,k),(\beta,j,\ell)\big)x_j^\beta(\ell)x_i^\alpha(k)
\, = \, 2\delta_{(\alpha,i,k),(\beta,j,\ell)}
\end{align}
for $(\alpha,i,k),(\beta,j,\ell) \in \Upsilon_m$. We endow $\mathcal
A_{\varepsilon_n}$ with the antilinear involution such that
$x_{i}^{\alpha}(k)^*=x_{i}^{\alpha}(k)$ for every $(\alpha,i,k)\in
\Upsilon_m$. If we equip $\Upsilon_m$ with the lexicographical order,
a basis of the linear space $\mathcal A_{\varepsilon_m}$ is given by
$x_\emptyset^{\eps_m } = \mathbf{1}_{\mathcal A_{\varepsilon_m}}$ and
the set of reduced words written in increasing order. Namely, elements
of the form
$$x_A^{\eps_m} \, = \, x_{i_1}^{\alpha_1}(k_1)\cdots
x_{i_s}^{\alpha_s}(k_s),$$ where $A = \{(\alpha_1,i_1,k_1), \ldots,
(\alpha_s,i_s,k_s)\} \subset \Upsilon_m$ with
$(\alpha_j,i_j,k_j)<(\alpha_{j+1},i_{j+1},k_{j+1})$ for $1 \le j \le
s-1$. For any such element we set $|x_A^{\eps_m}|=|A|=s$. Define the
tracial state on $\mathcal A_{\varepsilon_m}$ given by
$\tau_{\eps_m}(x_A^{\eps_m})=\delta_{\emptyset, A}$ for $A \subset
\Upsilon_m$. The given basis turns out to be orthonormal with respect
to the inner product $\langle x,y \rangle = \tau_{\eps_m}(x^*y)$. Let
$\mathcal A_{\varepsilon_m}$ act by left multiplication on the Hilbert
space $\mathcal H_{\mathcal A_{\varepsilon_m}} = (\mathcal
A_{\varepsilon_m}, \langle \cdot,\cdot\rangle)$ to get a faithful
$*$-representation of $\mathcal A_{\varepsilon_m}$ on $\mathcal
H_{\mathcal A_{\varepsilon_m}}$. We may endow $\mathcal
A_{\varepsilon_m}$ with the von Neumann algebra structure induced by
this representation and denote by $L_p(\mathcal A_{\varepsilon_m},
\tau_{\eps_m})$ the associated non-commutative $L_p$-space. At this
point, it is natural to define the
\emph{$\varepsilon_m$-Ornstein-Uhlenbeck semigroup} on $\mathcal
A_{\varepsilon_m}$ by
\begin{equation}\label{Ptepsn}
\mathcal{S}_{\eps_m,t}(x_A^{\eps_m}) \, = \, e^{-t|x_A^{\eps_m}|}
x_A^{\eps_m}.
\end{equation}
Biane extended hypercontractivity for fermions to these spin algebras
in \cite{Biane}
\begin{equation} \label{Biane's theorem matrix model}
\big\| \mathcal{S}_{\varepsilon_m,t}: L_p(\mathcal A_{\varepsilon_m})
\rightarrow L_q(\mathcal A_{\varepsilon_m}) \big\| \, = \, 1
\ \Leftrightarrow \ t \, \ge \, \frac12 \log \frac{q-1}{p-1},
\end{equation}
whenever $1 < p \le q < \infty$. We will also use the following direct
consequence of Biane's result. Namely, given $1 \le p < \infty$ and $r
\in \Z_+$ we may find constants $C_{p,r} > 0$ such that the following
inequality holds uniformly for all $m \geq 1$ and all homogeneous
polynomials $P$ of degree $r$ in $|\Upsilon_m|$ noncommutative
indeterminates satisfying \eqref{relations free} and written in
reduced form
\begin{equation} \label{Biane 2-p}
\Big\|P\big((x_i^\alpha(k))_{(\alpha,i,k) \in \Upsilon_m} \big)
\Big\|_{L_p(\mathcal A_{\varepsilon_m})} \, \le \, C_{p,r} \Big\| P
\big((x_i^\alpha(k))_{(\alpha,i,k) \in \Upsilon_m} \big)
\Big\|_{L_2(\mathcal A_{\varepsilon_m})}.
\end{equation}
According to \eqref{Biane's theorem matrix model}, it is
straightforward to show that we can take $C_{p,r}=(p-1)^{r/2}$.

\subsection{A central limit theorem}

In order to approximate the free product $\M$ of Clifford algebras, we
need to choose the commutation/anticommutation relations
randomly. More precisely, we consider a probability space $(\Omega,
\mu)$ and a family of independent random variables $$\varepsilon
\big((\alpha,i,k),(\beta,j,\ell) \big): \Omega \rightarrow \{-1,1\}
\quad \mbox{for} \quad (\alpha,i,k)<(\beta,j,\ell)$$ which are
distributed as follows
\begin{align}\label{sign choice}
\mu \Big(\varepsilon\big(( \alpha,i,k),(\beta,j,\ell)\big)=-1\Big) \,
= \, \begin{cases} 1 & \mbox{if } \alpha=\beta, \\ 1/2 & \text{if }
  \alpha\neq\beta. \end{cases}
\end{align}
In particular, this means that all the generators
$(x_i^\alpha(k))_{i\in [d],k\in[m]}$ anticommute for $\alpha \in [n]$
fixed and all $m \ge 1$. Therefore, the algebra $\A_{\eps_m}^\alpha$
generated by them is isomorphic to $\mathcal C(\R^{dm})$. Formally, we
have a matrix model for each $\omega \in \Omega$. In this sense, the
generators $x_i^\alpha(k)$ and the algebras $\A_{\eps_m}^\alpha$ are
also functions of $\omega$. In order to simplify the notation, we will
not specify this dependence unless it is necessary for clarity in the
exposition. Define also the algebra $$\tilde{\mathcal
  A}^\alpha_{\varepsilon_m} \, = \, \Big\langle \tilde{x}_i^\alpha(m)
\, \big| \ i \in [d] \Big\rangle$$ with generators given
by $$\tilde{x}_i^\alpha(m) \, = \, \frac{1}{\sqrt{m}}\sum_{k=1}^m
x_i^\alpha(k).$$

\begin{lemma}\label{single algebra}
The von Neumann algebra $\tilde{\mathcal A}^\alpha_{\varepsilon_m}$ is
canonically isomorphic to $\mathcal C(\R^d)$.
\end{lemma}

\dem It suffices to prove that the generators verify the CAR
relations. All of them are self-adjoint since the same holds for the
$x_i^\alpha$'s. Since $\alpha$ is fixed, our choice (\ref{sign
  choice}) of the sign function $\varepsilon$ and \eqref{relations
  free} give \\ \vskip-20pt \null${}$ \hfill $\displaystyle
\tilde{x}_i^\alpha(m)
\tilde{x}_j^\alpha(m)+\tilde{x}_j^\alpha(m)\tilde{x}_i^\alpha(m) \, =
\, \frac{1}{m}\sum_{k=1}^m\sum_{\ell=1}^m
x_i^\alpha(k)x_j^\alpha(\ell)+x_j^\alpha(\ell)x_i^\alpha(k) \, = \,
2\delta_{ij}.$ \fin

We will denote by $\Pi(s)$ the set of all partitions of $[s] =
\{1,2,\ldots,s\}$. Given $\sigma,\pi \in \Pi(s)$, we will write
$\sigma\leq \pi$ if every block of the partition $\sigma$ is contained
in some block of $\pi$. We denote by $\sigma_0$ the smallest
partition, in which every block is a singleton. Given an $s$-tuple
$\ibar=(i_1,\cdots, i_s)\in [N]^s$ for some $N$, we can define the
partition $\sigma(\ibar)$ associated to $\ibar$ by imposing that two
elements $j,k \in [s]$ belong to the same block of $\sigma(\ibar)$ if
and only if $i_j=i_k$. We will denote by $\Pi_2(s)$ the set of all
pair partitions. That is, partitions $\sigma=\{V_1,\cdots, V_{s/2}\}$
such that $|V_j|=2$ for every block $V_j$. In this case, we will write
$V_j = \{e_j,z_j\}$ with $e_j< z_j$ so that $e_1< e_2 < \cdots <
e_{s/2}$. For a pair partition $\sigma\in \Pi_2(s)$ we define the set
of crossings of $\sigma$ by $$I(\sigma) \, = \, \Big\{ (k,\ell) \,
\big| \ 1 \le k, \ell \le s, \ e_k < e_\ell < z_k < z_\ell \Big\}.$$
Moreover, given an $s$-tuple $\underline{\alpha} = (\alpha_1,\ldots,
\alpha_s)$ such that $\sigma \leq \sigma(\underline{\alpha})$, we can
define the set of crossings of $\sigma$ with respect to
$\underline{\alpha}$ by $I_{\underline{\alpha}}(\sigma)=\{(k,\ell) \in
I(\sigma) : \alpha_{e_k} \neq \alpha_{e_\ell}\}$. This notation allows
us to describe the moments of reduced words in $\M$ with a simple
formula. Indeed, the following lemma arises from \cite[Lemma
  2]{Speicher} and a simple induction argument like the one used below
to prove identity \eqref{ByInduction}.
\begin{lemma}\label{moments free pdct fermions}
If $\ibar \in [d]^s$ and $\underline{\alpha} \in [n]^s$ we
have $$\tau\big(x_{i_1}^{\alpha_1} \cdots x_{i_s}^{\alpha_s}\big) \, =
\, \delta_{s \in 2 \Z} \sum_{\substack{\sigma\in \Pi_2(s)\\ \sigma
    \leq \sigma(\ibar), \sigma(\underline{\alpha})
    \\ I_{\underline{\alpha}}(\sigma)=\emptyset}}(-1)^{\#I(\sigma)}.$$
\end{lemma}
We can now prove that the moments of the free product von Neumann
algebra $\M$ are the almost everywhere limit of the moments of our
matrix model. More explicitly, we find the following central limit
type theorem.
\begin{theorem}\label{moments}
If $\ibar \in [d]^s$ and $\underline{\alpha} \in [n]^s$ we
have $$\lim_{m\to \infty} \tau_{\eps_m}
\Big(\tilde{x}_{i_1}^{\alpha_1}(m)(\omega) \cdots
\tilde{x}_{i_s}^{\alpha_s}(m)(\omega) \Big) \, = \, \tau
\big(x_{i_1}^{\alpha_1} \cdots x_{i_s}^{\alpha_s} \big) \quad
\mbox{a.e.}$$
\end{theorem}
%
\dem We will first prove that the convergence holds in expectation.
For $\omega\in \Omega$ fixed, by developing and splitting the sum
according to the distribution we obtain
\begin{eqnarray} \label{sum with partitions}
\lefteqn{\hskip-15pt \tau_{\eps_m} \Big(
  \tilde{x}_{i_1}^{\alpha_1}(m)(\omega) \cdots
  \tilde{x}_{i_s}^{\alpha_s}(m)(\omega) \Big)} \\ [3pt] \nonumber & =
& \frac{1}{m^{s/2}} \sum_{\k \in [m]^s} \tau_{\eps_m}
\big(x_{i_1}^{\alpha_1}(k_1)(\omega) \cdots
x_{i_s}^{\alpha_m}(k_s)(\omega) \big) \\ [3pt] \nonumber &= &
\frac{1}{m^{s/2}} \sum_{\sigma \in \Pi(s)}
\underbrace{\sum_{\substack{\k \in [m]^s \\ \sigma (\k)=\sigma}}
  \tau_{\eps_m} \big( x_{i_1}^{\alpha_1}(k_1)(\omega)\cdots
  x_{i_s}^{\alpha_s}(k_s)(\omega)\big)}_{\mu_{\sigma}(\omega)}.
\end{eqnarray}
We claim that $$\lim_{m \to \infty} \frac{1}{m^{s/2}}
\mu_{\sigma}(\omega) \, = \, 0$$ for every $\sigma \in \Pi(s)
\backslash \Pi_2(s)$ and all $\omega \in \Omega$. Indeed, the upper
bound $\mu_{\sigma}(\omega) \le m^r$ holds when $\sigma$ has $r$
blocks since $|\tau_{\eps_m} ( x_{i_1}^{\alpha_1}(k_1)(\omega)\cdots
x_{i_s}^{\alpha_s}(k_s)(\omega))| \le 1$. Hence, the limit above
vanishes for $r < s/2$. It then suffices to show that the same limit
vanishes when $\sigma$ contains a singleton $\{j_0\}$. However, in
this case we have $$\tau_{\eps_m} \big(x_{i_1}^{\alpha_1}(k_1)(\omega)
\cdots x_{i_s}^{\alpha_s}(k_s)(\omega) \big) \, = \, 0$$ whenever
$\sigma (\k)=\sigma$ since the $j_0$-th term can not be
cancelled. This proves our claim. Hence, the only partitions which may
contribute in the sum (\ref{sum with partitions}) are pair partitions
$\sigma=\{\{e_1,z_1\}, \ldots ,\{e_\frac{s}{2},z_\frac{s}{2}\}\}$.  In
particular, if $s$ is odd we immediately obtain that the trace
converges to zero in (\ref{sum with partitions}).  Note that given
such a pair partition $\sigma$, we must have that $\sigma \leq
\sigma(\underline{\alpha})$ and $\sigma \leq \sigma(\ibar)$. Indeed,
if this is not the case we will have $i_{e_j}\neq i_{z_j}$ or
$\alpha_{e_j}\neq \alpha_{z_j}$ for some $j=1,2,\ldots, s/2$. Now, for
every $\k \in [m]^s$ such that $\sigma(\k)= \sigma$ we have
$k_{e_j}=k_{z_j}\neq k_\ell$ for every $\ell \neq e_j,z_j$. Thus, the
only way for the
elements $$x_{i_{e_j}}^{\alpha_{e_j}}(k_{e_j})(\omega) \quad
\mbox{and} \quad x_{i_{z_j}}^{\alpha_{z_j}}(k_{z_j})(\omega)$$ to
cancel is to match each other. Thus, we can assume that
$(\alpha_{e_j},i_{e_j})=(\alpha_{z_j},i_{z_j})$.

We have seen that the letters of our word should match in pairs. We
are now reduced to study the sign which arises from the
commutation/anticommutation relations to cancel all elements. Assume
that $\sigma$ has a crossing with respect to
$\underline{\alpha}=(\alpha_1,\ldots, \alpha_s)$. That is, there
exists $(k,\ell) \in I(\sigma)$ such that $\alpha_{e_k} \neq
\alpha_{e_\ell}$. Then we find that $$\mathbb{E}_\omega \tau_{\eps_m}
\big(x_{i_1}^{\alpha_1}(k_1)(\omega) \cdots
x_{i_s}^{\alpha_s}(k_s)(\omega) \big) \, = \, 0$$ for every
$(k_1,\ldots ,k_s)$ such that $\sigma(k_1,\ldots ,k_s) =
\sigma$. Indeed, define the sign
function $$\varepsilon^{\underline{\alpha}}_{(k,\ell)}:=\varepsilon\big((\alpha_{e_\ell},i_{e_\ell},k_{e_\ell}),(\alpha_{z_k},i_{z_k},k_{z_k})
\big).$$ If $\sigma$ has such a crossing, we obtain (among others)
this sign only once when canceling the letters associated to
$(\alpha_{e_k}, i_{e_k}, k_{e_k})$ and $(\alpha_{z_k}, i_{z_k},
k_{z_k})$ as well as $(\alpha_{e_\ell}, i_{e_\ell}, k_{e_\ell})$ and
$(\alpha_{z_\ell},i_{z_\ell}, k_{z_\ell})$. Furthermore, by
independence and since
$\mathbb{E}_\omega\varepsilon^{\underline{\alpha}}_{(k,\ell)}=0$ we
get
\begin{eqnarray*}
\lefteqn{\hskip-20pt \mathbb{E}_\omega \tau_{\eps_m}
  \big(x_{i_1}^{\alpha_1}(k_1)(\omega)\cdots
  x_{i_s}^{\alpha_s}(k_1)(\omega)\big)} \\ & = & \pm \,
\mathbb{E}_\omega \Big( \prod_{(k,\ell) \in
  I_{\underline{\alpha}}(\sigma)}
\varepsilon^{\underline{\alpha}}_{(k,\ell)} \Big) \ = \ \pm
\prod_{(k,\ell) \in I_{\underline{\alpha}}(\sigma)} \mathbb{E}_\omega
\varepsilon^{\underline{\alpha}}_{(k,\ell)}=0,
\end{eqnarray*}
where $\pm$ denotes a possible change of signs depending on the
crossings of $\sigma$. Then, we can also rule out these kind of
partitions and we can assume that $\sigma \in \Pi_2(s)$ is such that
$\sigma \le \sigma(\ibar), \sigma(\underline{\alpha})$ and
$I_{\underline{\alpha}}(\sigma)=\emptyset$. In this case, we do not
need to commute two letters $(\alpha,i,k)$ and $(\beta,j,\ell)$ with
$\alpha \neq \beta$. Hence we will obtain deterministic signs coming
from the commutations, which only depend on the number of crossings of
$\sigma$. More precisely, given $\sigma \in \Pi_2(s)$ satisfying the
properties above and $\k \in [m]^s$ such that $\sigma(\k) = \sigma$ we
have
\begin{equation} \label{ByInduction}
\tau_{\eps_m} \Big(x_{i_1}^{\alpha_1}(k_1)(\omega)\cdots
x_{i_s}^{\alpha_s}(k_s)(\omega) \Big)=(-1)^{\#I(\sigma)} \quad
\mbox{for every} \quad \omega.
\end{equation}
Indeed, this can be proved inductively as follows. Using that
$I_{\underline{\alpha}}(\sigma) = \emptyset$, there must exists a
connected block of consecutive numbers in $[s]$ so that the following
properties hold
\begin{itemize}
\item The letters in that block are related to a fixed $\alpha$.

\item The product of the letters in that block equals $\pm
  \mathbf{1}$.

\item The block itself is a union of pairs of the partition $\sigma
  \in \Pi_2(s)$.
\end{itemize}
If $\pi$ denotes the restriction of $\sigma$ to our distinguished
block ---well defined by the third property--- the sign given by the
second property equals $(-1)^{\#I(\pi)}$. After canceling this block
of letters, we may start again by noticing that
$I_{\underline{\beta}}(\sigma \setminus \pi) = \emptyset$ where
$\underline{\beta}$ is the restriction of $\underline{\alpha}$ to the
complement of our distinguished block. This allows to restart the
process. In the end we obtain $(-1)^{\#I(\sigma)}$ as desired. We
deduce that
\begin{eqnarray*}
\lefteqn{\lim_{m \to \infty} \mathbb{E}_\omega \tau_{\eps_m} \Big(
  \tilde{x}_{i_1}^{\alpha_1}(m)(\omega) \cdots
  \tilde{x}_{i_s}^{\alpha_s}(m)(\omega) \Big)} \\ [3pt] & = & \lim_{m
  \to \infty} \frac{1}{m^{s/2}} \mathbb{E}_\omega
\sum_{\substack{\sigma\in \Pi_2(s) \\ \sigma \le
    \sigma(\ibar),\sigma(\underline{\alpha})
    \\ I_{\underline{\alpha}}(\sigma)=\emptyset}} \sum_{\substack{\k
    \in [m]^s \\ \sigma(\k)=\sigma}} (-1)^{\#I(\sigma)} \ =
\ \sum_{\substack{\sigma\in \Pi_2(s) \\ \sigma \le \sigma(\ibar),
    \sigma\le \sigma(\underline{\alpha})
    \\ I_{\underline{\alpha}}(\sigma)=\emptyset}} (-1)^{\#I(\sigma)}.
\end{eqnarray*}
Here we have used that $$\lim_{m \to \infty} \frac{|\{\k \in [m]^s :
  \sigma(\k) = \sigma\}|}{m^{s/2}} \, = \, \lim_{m \to \infty}
\frac{m(m-1)\cdots (m-\frac{s}{2}+1)}{m^{s/2}} \, = \, 1.$$ By Lemma
\ref{moments free pdct fermions}, this proves convergence in
expectation and completes the first part of the proof. It remains to
prove almost everywhere convergence in $\omega$.  Let us define the
random variables $$X_m(\omega) \, = \, \tau_{\eps_m} \Big(
\tilde{x}_{i_1}^{\alpha_1}(m) \cdots \tilde{x}_{i_s}^{\alpha_s}(m)
\Big).$$ By the dominated convergence theorem, it suffices to
show $$\lim_{m \rightarrow \infty} \mu \Big( \big\{\sup_{M\geq
  m}\big|X_M-\mathbb{E}_\omega[X_M]\big|\geq \alpha\big\} \Big) \, =
\, 0$$ for every $\alpha> 0$. According to Tchebychev's inequality, we
find $$\mu\Big(\big\{\sup_{M\geq
  m}\big|X_M-\mathbb{E}_\omega[X_M]\big|\geq \alpha\big\}\Big) \, \le
\, \frac{1}{\alpha^2}\sum_{M=m}^\infty V[X_M],$$ where
$V[X_M]=\mathbb{E}_\omega[X_M^2]- (\mathbb{E}_\omega[X_M] )^2$ denotes
the variance of $X_M$. We will prove the upper bound $V[X_M] \leq
C(s)/M^2$ for every $M$, for some contant $C(s)$ depending only on the
length $s$. This will suffice to conclude the argument. To this end we
write
\begin{equation} \label{variance}
V[X_M] \, = \, \frac{1}{M^s} \sum_{\substack{\sigma, \pi \in \Pi(s)}}
\, \sum_{\substack{\k \, : \, \sigma(\k) = \sigma \\ \l \, : \, \sigma
    \hskip1pt (\l) \hskip1pt = \pi}} D_{\k,\l},
\end{equation}
where
\begin{eqnarray*}
D_{\k,\l} & = & \mathbb{E}_\omega \Big[ \tau_{\eps_m} \big(
  x_{i_1}^{\alpha_1}(k_1)(\omega) \cdots
  x_{i_s}^{\alpha_s}(k_s)(\omega)\big) \tau_{\eps_m} \big(
  x_{i_1}^{\alpha_1}(\ell_1)(\omega) \cdots
  x_{i_s}^{\alpha_s}(\ell_s)(\omega) \big) \Big] \\ & - &
\mathbb{E}_\omega \Big[ \tau_{\eps_m}
  \big(x_{i_1}^{\alpha_1}(k_1)(\omega) \cdots
  x_{i_s}^{\alpha_s}(k_s)(\omega) \big) \Big] \mathbb{E}_\omega \Big[
  \tau_{\eps_m} \big(x_{i_1}^{\alpha_1}(\ell_1)(\omega) \cdots
  x_{i_s}^{\alpha_s}(\ell_s)(\omega) \big) \Big]
\end{eqnarray*}
for $\k = (k_1,\ldots, k_s)$ and $\l=(\ell_1,\ldots, \ell_s)$. Now,
reasoning as above one can see that whenever $\sigma$ or $\pi$ has a
singleton, all the corresponding terms in the sum \eqref{variance} are
equal to zero. Thus, we may write $\sigma=\{V_1,\ldots,
V_{r_\sigma}\}$ and $\pi=\{W_1,\ldots, W_{r_\pi}\}$ with
$r_\sigma,r_\pi \le \frac{s}{2}$. If neither $\sigma$ nor $\pi$ are
pair partitions, we will have $r_\sigma,r_\pi \le \frac{s}{2}-1$ and
the part of the sum in \eqref{variance} corresponding to these pairs
$(\sigma,\pi)$ can be bounded above in absolute value by $C(s) / M^2$
as desired. Then, it remains to control the rest of the terms in
\eqref{variance}. To this end, we assume that $\sigma$ is a pair
partition. Actually, a cardinality argument as before allows us to
conclude that $\pi$ must be either a pair partition or a partition
with all blocks formed by two elements up to a possible four element
block. In the following, we will explain how to deal with the case in
which $\pi$ is a pair partition. The other case can be treated exactly
in the same way, being actually even easier by cardinality
reasons. Let us fix two pair partitions $\sigma$ and $\pi$ and let us
consider $\k = (k_1,\ldots, k_s)$ and $\l = (\ell_1,\ldots, \ell_s)$
such that $\sigma(\k)=\sigma$ and $\sigma(\l)=\pi$. When rearranging
the letters in the traces defining $D_{\k,\l}$, the deterministic
signs ---$\alpha = \beta$ in \eqref{sign choice}--- do not have any
effect in the absolute value of $D_{\k,\l}$. On the other hand, the
random signs ---$\alpha \neq \beta$ in \eqref{sign choice}--- makes
the second term of $D_{\k,\l}$ vanish. Thus, $D_{\k,\l} \neq 0$ if and
only if $I_{\underline{\alpha}}(\sigma) \neq \emptyset \neq
I_{\underline{\alpha}}(\pi)$ and we obtain the same random signs
coming from crossings in $I_{\underline{\alpha}}(\sigma)$ and
$I_{\underline{\alpha}}(\pi)$. In particular, we should find at least
two signs $$\eps \big( (\alpha_p, i_p, k_p),
(\alpha_q,i_q,k_q)\big)(\omega) \quad (\alpha_p \hskip0.5pt \neq
\hskip0.5pt \alpha_q) \quad \mbox{from} \quad
x_{i_1}^{\alpha_1}(k_1)(\omega)\cdots
x_{i_s}^{\alpha_s}(k_s)(\omega),$$ $$\eps \big((\alpha_u,i_u,\ell_u),
(\alpha_v,i_v,\ell_v)\big) (\omega) \quad (\alpha_u \neq \alpha_v)
\quad \mbox{from} \quad x_{i_1}^{\alpha_1}(\ell_1) \hskip0.5pt
(\omega)\cdots x_{i_s}^{\alpha_s}(\ell_s) \hskip0.5pt (\omega).$$ By
independence, this implies that $$\big\{ (\alpha_p, i_p, k_p),
(\alpha_q,i_q,k_q) \big\} \, = \, \big\{ (\alpha_u, i_u, \ell_u),
(\alpha_v,i_v,\ell_v) \big\}.$$ Moreover, since we also need $\sigma
\leq \sigma(\underline{\alpha})$ for non-vanishing terms, we can
conclude that $k_p \neq k_q$ and $\ell_u \neq \ell_v$. Therefore, the
sets $\{k_1,\ldots, k_s\}$ and $\{\ell_1,\ldots, \ell_s\}$ must have
four elements (corresponding to two different blocks) in common. This
implies that the part of the sum in \eqref{variance} corresponding to
pairs $(\sigma,\pi)$ of pair partitions is bounded above by $$C'(s)
\frac{M^{s/2} M^{(s-4)/2}}{M^s} \, = \, \frac{C'(s)}{M^2}$$ for a
certain constant $C(s)'$ as we wanted. This completes the proof. \fin

Let $x$ be a word in the reduced free product of Clifford algebras
$\M$, which written in reduced form is given by \eqref{general element
  free productI}. In what follows, we will associate to $x$ an element
$\tilde{x}(m)$ in $\mathcal A_{\varepsilon_m}$ given by
\begin{equation}\label{general element free productI with n}
\tilde{x}(m) \, = \, \tilde{x}_{A_1}^{\alpha_1}(m) \cdots
\tilde{x}_{A_\ell}^{\alpha_\ell}(m).
\end{equation}
If we develop $x$ as in \eqref{general element free productII}, then
we can write $\tilde{x}(m)$
as $$\overbrace{\tilde{x}_{i_1}^{\alpha_1}(m) \cdots
  \tilde{x}_{i_{s_1}}^{\alpha_1}(m)}^{\tilde{x}_{A_1}^{\alpha_1}(m)}\overbrace{\tilde{x}_{i_{s_1+1}}^{\alpha_2}(m)
  \cdots
  \tilde{x}_{i_{s_1+s_2}}^{\alpha_2}(m)}^{\tilde{x}_{A_2}^{\alpha_2}}
\cdots \overbrace{\tilde{x}_{i_{s_1+\cdots
      +s_{\ell-1}+1}}^{\alpha_\ell}(m) \cdots \tilde{x}_{i_{s_1+\cdots
      +s_\ell}}^{\alpha_\ell}(m)}^{\tilde{x}_{A_\ell}^{\alpha_\ell}(m)}.$$

\subsection{Hypercontractivity bounds}

In this subsection we prove Theorem B. The result below can be obtained
following verbatim the proof of \cite[Lemma 4]{Biane} just replacing
Theorem 7 there by Theorem \ref{moments} above.

\begin{lemma}\label{approx norm}
If $p\geq 1$, we have $$\lim_{m \to \infty} \Big\| \summ_j \rho_j
\tilde{x}_j(m) \Big\|_{L_p(\mathcal A_{\varepsilon_m})} \, = \, \Big\|
\summ_j \rho_j x_j \Big\|_{L_p(\M)} \quad \mbox{a.e.}$$ for any
finite linear combination $\sum_j \rho_j x_j$ of reduced words in
the free product $\M$.
\end{lemma}

\begin{lemma}\label{key lemma}
Given $x$ a reduced word in the free product $\M$, let $\tilde{x}(m)$
be the element in $\mathcal A_{\varepsilon_m}$ associated to $x$ as in
\eqref{general element free productI with n}. Then, there exists a
decomposition $\tilde{x}(m) = \tilde{x}_1(m) + \tilde{x}_2(m)$ with
the following properties
\begin{itemize}
\item[i)] $\langle \tilde{x}_1(m), \tilde{x}_2(m)\rangle=0$ a.e.,

\vskip3pt

\item[ii)]
  $\mathcal{S}_{\varepsilon_n,t}(\tilde{x}_1(m))=e^{-t|x|}\tilde{x}_1(m)$,

\vskip3pt

\item[iii)] $\displaystyle\lim_{m\rightarrow
  \infty}\|\tilde{x}_1(m)\|_{L_2(\mathcal A_{\varepsilon_m})}=1$ a.e.
\end{itemize}
In particular, we deduce that $$\lim_{m\rightarrow \infty}
\|\tilde{x}_2(m)\|_{L_2(\mathcal A_{\varepsilon_m})} \, = \, 0 \quad
\mbox{a.e}.$$
\end{lemma}

\dem If we set $s = |x|$ and $\sigma_0$ denotes the singleton
partition, define
\begin{eqnarray*}
\tilde{x}_1(m)(\omega) & = & \frac{1}{m^{s/2}} \sum_{\substack{\k \in
    [m]^s \\ \sigma(\k) = \sigma_0}} x_{i_1}^{\alpha_1}(k_1)(\omega)
\cdots x_{i_s}^{\alpha_s}(k_s)(\omega), \\ [5pt]
\tilde{x}_2(m)(\omega) & = & \frac{1}{m^{s/2}} \sum_{\sigma \in \Pi(s)
  \setminus \{\sigma_0\}} \sum_{\substack{\k \in[m]^s \\ \sigma(\k) =
    \sigma}} x_{i_1}^{\alpha_1}(k_1)(\omega) \cdots
x_{i_s}^{\alpha_s}(k_s)(\omega).
\end{eqnarray*}
Clearly $\tilde{x}(m) = \tilde{x}_1(m) + \tilde{x}_2(m)$ point wise
and $\|\tilde{x}(m)\|_{L_2(\mathcal A_{\varepsilon_m})}=1$. Property
i) is easily checked. Indeed, consider $\k,\underline{\ell}\in [m]^s$
with $\sigma(\k)=\sigma_0$ and $\sigma(\underline{\ell}) \in \Pi(s)
\setminus \{\sigma_0\}$. Since the $k_i$'s are all distinct and the
$\ell_i$'s are not we must have $$\tau_{\eps_m}
\Big(x_{i_1}^{\alpha_1}(k_1)(\omega)\cdots
x_{i_s}^{\alpha_s}(k_s)(\omega)x_{i_s}^{\alpha_s}(\ell_s)(\omega)\cdots
x_{i_1}^{\alpha_1}(\ell_1)(\omega)\Big) \, = \, 0.$$ The second
property comes from the definition of the semigroup \eqref{Ptepsn} and
the fact that for every $\k$ with $\sigma(\k)=\sigma_0$, we have no
cancellations. Now it remains to show that $$\lim_{m \rightarrow
  \infty} \frac{1}{m^s} \sum_{\substack{\k, \l \in [m]^s
    \\ \sigma(\k)=\sigma(\l)=\sigma_0}} \tau_{\eps_m}
\Big(\underbrace{x_{i_1}^{\alpha_1}(k_1)\cdots
  x_{i_s}^{\alpha_s}(k_s)}_{x_{\ibar}^{\underline{\alpha}}(\k)}
\underbrace{x_{i_s}^{\alpha_s}(\ell_s)\cdots
  x_{i_1}^{\alpha_1}(\ell_1)}_{x_{\ibar}^{\underline{\alpha}}(\l)^*}\Big)
\, = \, 1.$$ Indeed, if $\{k_1,\ldots, k_s\}\neq\{\ell_1,\ldots,
\ell_s\}$ the trace clearly vanishes and it suffices to consider the
case $\{k_1,\ldots, k_s\}=\{\ell_1,\ldots, \ell_s\}$. Note that, the
trace above is different from $0$ if and only if $(\alpha_j,i_j,k_j) =
(\alpha_{\beta(j)}, i_{\beta(j)}, \ell_{\beta(j)})$ for some
permutation $\beta \in S_s$ and every $1 \le j \le s$. If we assume
$k_s\neq \ell_s$, we get $(\alpha_j, i_j, k_j) = (\alpha_s, i_s,
\ell_s)$ for certain $j < s$. This means that
$x_{i_j}^{\alpha_j}(k_j)$ and $x_{i_s}^{\alpha_s}(k_s)$ belong to
different $\alpha$-blocks since the $i_j$'s are pairwise distinct in a
fixed $\alpha$-block. Thus, to cancel these elements we must cross a
$\beta$-block with $\beta \neq \alpha_s$. Since the $k$'s are all
different, the $\eps$-signs corresponding to these commutations appear
just once. We can argue in the same way for every $1\leq j\leq s$ and
conclude that
$$\mathbb{E}_\omega \, \tau_{\eps_m} \Big( x_{i_1}^{\alpha_1}(k_1)
\cdots x_{i_s}^{\alpha_s}(k_s) x_{i_s}^{\alpha_s}(\ell_s) \cdots
x_{i_1}^{\alpha_1}(\ell_1) \Big) \, = \, 0$$ unless $k_j = \ell_j$ for
all $1\leq j\leq s$. Therefore
\begin{eqnarray*}
\lim_{m \rightarrow \infty}\mathbb{E}_\omega
\|\tilde{x}_1(m)\|_{L_2(\mathcal A_{\varepsilon_m})}^2 & = & \lim_{m
  \rightarrow \infty} \frac{1}{m^s} \sum_{\substack{\k \in [m]^s
    \\ k_i\neq k_j}} 1 \\ & = & \frac{m(m-1)\cdots (m-s+1)}{m^s} \ =
\ 1.
\end{eqnarray*}
Finally, arguing as in the proof of Theorem \ref{moments} we see that
the same limit holds for almost every $\omega\in \Omega$.  This proves
iii). The last assertion follows from i), iii) and the identity
$\|\tilde{x}(m)\|_2=1$. The proof is complete. \fin

\begin{lemma}\label{approx normII}
If $p\geq 1$, we have $$\lim_{m \to \infty} \Big\|
\mathcal{S}_{\eps_m,t} \Big( \summ_j \rho_j \tilde{x}_j(m) \Big)
\Big\|_{L_p(\mathcal A_{\varepsilon_m})} \, = \, \Big\|
\mathcal{O}_{\M,t} \Big( \summ_j \rho_j x_j \Big) \Big\|_{L_p(\M)}
\quad \mbox{a.e.}$$ for any finite linear combination $\sum_j
\rho_j x_j$ of reduced words in the free product $\M$.
\end{lemma}

\dem According to Lemma \ref{key lemma}, we have $$\lim_{m \rightarrow
  \infty} \Big\|\Big(\mathcal{S}_{\varepsilon_m,t}- e^{-t|x|}
\mathbf{1}_{\mathcal
  A_{\varepsilon_m}}\Big)\big(\tilde{x}(m)\big)\Big\|_{L_2(\mathcal
  A_{\varepsilon_m})} \, = \, 0 \quad \mbox{a.e.}$$ for any reduced
word $x\in \M$ and the associated $\tilde{x}(m)$'s $\in \mathcal
A_{\varepsilon_m}$ given by \eqref{general element free productI with
  n}. Thus $$\lim_{m \rightarrow \infty}
\Big\|\mathcal{S}_{\varepsilon_m,t} \Big(\summ_j \rho_j
\tilde{x}_j(m) \Big) - \summ_j e^{-t|x_j|} \rho_j \tilde{x}_j(m)
\Big\|_{L_2(\mathcal A_{\varepsilon_m})} \, = \, 0 \quad \mbox{a.e.}$$
Then \eqref{Biane 2-p} implies that the same limit vanishes in the
norm of $L_p(\mathcal{A}_{\eps_m})$. On the other hand, since
$\mathcal{O}_{\M,t} (x_j)=e^{-t|x_j|}x_j$, the assertion follows from
Lemma \ref{approx norm}. \fin

\demB Let $1<p\leq q <\infty$.  By construction, the algebraic free
product $A$ is a weak-$*$ dense involutive subalgebra of $\M$. In
particular, it is dense in $L_p(\M)$ for every $p<\infty$. Given a
finite sum $z = \sum_j \rho_j x_j \in A$, consider the
corresponding sum $\tilde{z}(m) = \sum_j \rho_j
\tilde{x}_j(m)\in\mathcal A_{\varepsilon_m}$ following \eqref{general
  element free productI with n}. Given any $t \ge \frac12 \log
(q-1/p-1)$, we may apply Lemmas \ref{approx norm} and \ref{approx
  normII} in conjunction with Biane's theorem \eqref{Biane's theorem
  matrix model} to conclude
\begin{eqnarray*}
\|\mathcal{O}_{\M,t}(z)\|_{L_q(\M)} & = & \lim_{m \to \infty} \big\|
\mathcal{S}_{\eps_m,t}(\tilde{z}(m))
\big\|_{L_q(\mathcal{A}_{\eps_m})} \\ & \le & \lim_{m \to \infty} \|
\tilde{z}(m)) \|_{L_p(\mathcal{A}_{\eps_m})} = \|z\|_{L_p(\M)}.
\end{eqnarray*}
The necessity of the condition $t \ge \frac12 \log(q-1/p-1)$ was
justified above. \fin

\subsection{Further comments} \label{remark generalization}

Note that the argument we have used in the proof of Theorem B still
works in a more general setting. More precisely, we may replace the
fermion algebras $\M_\alpha=\mathcal C(\R^d)$ by spin system algebras
$\mathcal{A}_\alpha$, where the generators $x_i^\alpha$ satisfy
certain commutation and anticommutation relations given by a sign
$\eps^\alpha$ as follows $$x_i^\alpha
x_j^\alpha-\eps^\alpha(i,j)x_j^\alpha x_i^\alpha=2\delta_{ij} \quad
\mbox{ for } 1\leq i,j\leq d.$$ Indeed, we just need to replace
\eqref{sign choice} by $$\mu \Big(\varepsilon \big((
\alpha,i,k),(\beta,j,\ell)\big)=-1\Big)= \begin{cases}
  \eps^\alpha(i,j) & \mbox{if} \ \alpha=\beta, \\ 1/2 & \mbox{if}
  \ \alpha \neq \beta. \end{cases}$$ This yields optimal time
hypercontractivity bounds for the Ornstein-Uhlenbeck semigroup on the
free product of spin matrix algebras. An additional application of
Speicher's central limit theorem allows us to obtain optimal
hypercontractivity estimates for the Ornstein-Uhlenbeck semigroup on
the free product of $q$-deformed algebras $\Gamma_q$, $-1\leq q \leq
1$.

\begin{remark}
\emph{Slight modifications in \eqref{sign choice} lead to
von Neumann algebras which are still poorly understood.
For instance, let us fix a function $f:[1,n]\times [1,n]\to
[-1,1]$ which is symmetric and assume that
 \[ \mu \big( \{\eps((\alpha,i,k),(\beta,j,\ell)) = +1\} \big)
 = \frac{1+f(\alpha,\beta)}{2} .\]
As usual we will assume that all the random variables
$\eps(x,y)$ are independent. Then it is convenient to
first calculate expectation of the joint moments of
 \[ \tilde{x}_i^{\alpha}(m) = \frac{1}{\sqrt{m}}
 \sum_{k=1}^m x_i^{\alpha}(k) .\]
Again, only the pair partitions survive and we get
\begin{eqnarray*}
\lim_{m \to \infty}\mathbb{E}_\omega \tau_{\varepsilon_m}(\tilde{x}_{i_1}^{\alpha_1}(m)
\cdots \tilde{x}_{i_s}^{\alpha_s}(m))   &=&
\sum_{\substack{\sigma \in \Pi_2(s)\\ \sigma \leq \sigma(\underline{i}), \sigma(\underline{\alpha})}}\prod_{(k,\ell)\in I(\sigma)}f(\alpha_{e_k},\alpha_{e_\ell}).
\end{eqnarray*}
As above, we will have hypercontractivity with the
optimal constant for the limit gaussian systems  (they indeed
produce a tracial von Neumann algebra). As an
illustration, let us consider $n=2$, $q_1,q_2 \in
[-1,1]$, $f(1,1)= q_1q_2$ and
$f(1,2)=f(2,1)=f(2,2)=q_2$. We deduce immediately that}
 \begin{itemize}

  \item[i)] \emph{The von Neumann subalgebra
      generated by $$x_i^1 = \lim_m
      \tilde{x}_i^1(m),$$ for $i=1, \ldots,\, d$ is
      isomorphic to $\Gamma_{q_1q_2}(\R^d)$, generated by $d$
      $q_1q_2$-gaussians.}

   \vskip3pt

    \item[ii)] \emph{The von Neumann subalgebra
        generated by $$x_i^2 = \lim_m
        \tilde{x}_i^2(m),$$ for $i=1, \ldots,\, d$ is
        isomorphic to $\Gamma_{q_2}(\R^d)$, generated by $d$
        $q_2$-gaussians.}

 \vskip3pt

  \item[iii)] \emph{Let $A\subset [s]$ and
let $y_i=x_{j_i}^1$ for $i\in A$ (and $\alpha_i=1$) and $y_i =
x_{j_i}^2$ ($\alpha_i=2$) otherwise.
 Let $\eta_0$ be the partition of $[s]$ defined by the possible
values of $(j_i,\alpha_i)$. Then we get
   \[ \tau(y_{1}y_{2}\cdots
   y_{s})
   = \sum_{\eta_0\geq \sigma \in \Pi_2(s)} q_1^{{\rm
   \scriptsize inversion}(\sigma |A)}   q_2^{{\rm
   \scriptsize inversion}(\sigma)}
   . \]
Here $\si|A$ is the restriction of $\si$ to $A$ where we count only inversions inside $A$. This construction is considered in \cite{JC} for constructing new Brownian motions.}
\end{itemize}
\emph{We see that we can combine different $q$ gaussian
random variables in one von Neumann algebra with a
prescribed interaction behaviour. With this method we
recover the  construction from \cite{JC}
of a non-stationary Brownian motion $B_t$. Indeed one
can choose $0=t_0<t_1<\cdots < t_d$ such that $B_t$ is an
abstract Brownian motion \cite{JC}
and the random variables $s_t(j)=B_t-B_{t_j}$ are
$q_0\cdots q_j$-Brownian motions. In this construction we needed a $q_1$-Brownian motion over a $q_2$-Brownian motion and hence the choice of the product $q_1q_2$ above. Although it is no
longer trivial to determine the number operator, we see
that hypercontractivity is compatible with
non-stationarity. The algebras generated for arbitrary
symmetric  $f$ could serve as models for $q_1$-products
over $q_2$-products, although in general there is no
$q$-product of arbitrary von Neumann algebras.}
\end{remark}

\section{{\bf The free Poisson semigroup}}\label{Section: proof of theorem A}

In this section we prove Theorem A and optimal hypercontractivity
for linear combinations of words in $\F_n$ with length lower than or equal to
1. Let us start with a trigonometric identity, which follows
from the binomial theorem and the identity $2\cos x = e^{ix} +
e^{-ix}$ $$(\cos x)^m \, = \, \frac{1}{2^{m-1}} \sum_{0 \le k \le
  [\frac{m}{2}]} {m \choose k}
\frac{\cos((m-2k)x)}{2^{\delta_{m,2k}}}.$$ Let $g_j$ denote one of the
generators of $\F_n$. Identifying $\lambda(g_j)$ with $\exp(2\pi i
\cdot)$, the von Neumann algebra generated by $\lambda(g_j)$ is
$\mathcal{L}(\Z)$ and the previous identity can be rephrased as
follows for $u_j=\lambda(g_j)$
\begin{equation} \label{Trig-Fn}
\big( u_j+u_j^* \big)^m \, = \, \sum_{0\leq k\leq [\frac{m}{2}]}{m
  \choose k} v_{j,m-2k},
\end{equation}
with $v_{j,k}=u_j^k+(u_j^*)^k$ for every $k \geq 1$ and
$v_0=\mathbf{1}$. We will also need a similar identity in
$\mathbb{G}_{2n}$. Let $z_1, z_2, \ldots, z_{2n}$ denote the canonical
generators of $\mathbb{G}_{2n}$, take $x_j=\lambda(z_j)$ for $1\leq j\leq 2n$
and consider the operators $a_{j,0} = \mathbf{1}, b_{j,0} = 0$ and
\begin{equation}\label{defab}
 a_{j,k} \, = \,
\underbrace{x_{2j-1} x_{2j} x_{2j-1}\cdots}_k \qquad
\mbox{,} \qquad b_{j,k} \, = \,\underbrace{x_{2j} x_{2j-1}
  x_{2j} \cdots}_k.
\end{equation}
If we set $\zeta_j = u_j + u_j^*$ and $\psi_j = x_{2j-1} + x_{2j}$, let us consider
the $*$-homomorphism $\Lambda: \mathcal{A}_{sym}^n \to
\mathcal{L}(\mathbb{G}_{2n})$ determined by $\Lambda(\zeta_j) =
\psi_j$. The result below can be proved by induction summing by
parts.

\begin{lemma}\label{lemma symmetric}
If $m \ge 0$, we find $$\big( x_{2j-1} + x_{2j} \big)^m \, = \,
\sum_{0 \leq k \leq [\frac{m}{2}]} {m \choose k} \big( a_{j,m-2k} +
b_{j,m-2k} \big).$$ Moreover, $v_{j,k} \in \langle u_j+u_j^* \rangle$
and we have $\Lambda(v_{j,k}) = a_{j,k} + b_{j,k}$ for every $k \geq
0$.
\end{lemma}

\demA As observed in the Introduction, the group von Neumann algebra
$\mathcal{L}(\Z_2)$ is $*$-isomorphic to the Clifford algebra
$\mathcal{C}(\R)$. Moreover, the Poisson and Ornstein-Uhlenbeck
semigroups coincide in this case. In particular, the first assertion
follows from $\mathcal{L}(\mathbb{G}_n)=\mathcal{L}(\Z_2) *\cdots *
\mathcal{L}(\Z_2) \simeq \mathcal C(\R) * \cdots * \mathcal C(\R)$, by
applying Theorem B with $d=1$. To prove the second assertion, we
consider the injective group homomorphism determined by $$\Phi: g_j
\in \F_n \, \mapsto \, x_{2j-1}x_{2j} \in \mathbb{G}_{2n}.$$ This map
clearly lifts to an isometry $L_p(\mathcal{L}(\F_n)) \to
L_p(\mathcal{L}(\mathbb{G}_{2n}))$ for all $p\geq 1$. Moreover, since
$|\Phi(g)|=2|g|$, we see that $\Phi$ intertwines the corresponding
free Poisson semigroup up to a constant $2$. More precisely, $\Phi
\circ \mathcal{P}_{\F_n,t} \, = \, \mathcal{P}_{\mathbb{G}_{2n}, t/2}
\circ \Phi$ for all $t > 0$. Hence, if $1<p\leq q <\infty$ and $f \in
L_p(\mathcal{L}(\F_n))$, we obtain from the result just proved
that $$\big\|\mathcal{P}_{\F_n,t} f \big\|_{L_q(\mathcal{L}(\F_n))} =
\big\|(\mathcal{P}_{\mathbb{G}_{2n},t/2} \circ \Phi) f
\big\|_{L_q(\mathcal{L}(\mathbb{G}_{2n}))} \le \|\Phi
f\|_{L_q(\mathcal{L}(\mathbb{G}_{2n})}=\|f\|_{L_q(\mathcal{L}(\F_n))},$$
whenever $t \ge \log (q-1/p-1)$. It remains to prove the last
assertion iii). The necessity of the condition $t \ge \frac12 \log
(q-1/p-1)$ can be justified following Weissler argument in \cite[pp
  220]{Weissler}. Therefore, we just need to prove
sufficiency. According to \cite{NiSp}, $\chi_{[-2,2]}(s)/\pi
\sqrt{4-s^2}$ is the common distribution of $\zeta_j$ and
$\psi_j$. Moreover, since both families of variables are free, the
tuples $(\zeta_1,\ldots, \zeta_n)$ and $(\psi_1,\ldots, \psi_n)$ must
have the same distribution too. Therefore, for every polynomial $P$ in
$n$ non-commutative variables we
have $$\big\|P(\zeta_1,\ldots,\zeta_n)\big\|_{L_p(\mathcal A_{sym}^n)}
\, = \,
\big\|P(\psi_1,\ldots,\psi_n)\big\|_{L_p(\mathcal{L}(\mathbb{G}_{2n}))}$$
for every $1\leq p\leq \infty$. In particular, the $*$-homomorphism
$\Lambda: \mathcal{A}_{sym}^n \to \mathcal{L}(\mathbb{G}_{2n})$
determined by $\Lambda(\zeta_j)=\psi_j$ for every $1\leq j\leq n$
extends to an $L_p$ isometry for every $1\leq p\leq \infty$. We claim
that $$\Lambda \big( \mathcal{P}_{\F_n,t}(P(\zeta_1,\ldots, \zeta_n))
\big) \, = \, \mathcal{P}_{\mathbb{G}_{2n},t} \big(P(\psi_1,\ldots,
\psi_n) \big)$$ for every polynomial $P$ in $n$ non-commutative
variables. It is clear that the last assertion iii) of Theorem A
follows from our claim above in conjunction with the first assertion
i), already proved. By freeness of the semigroups involved and the
fact that $\Lambda$ is a $*$-homomorphism, it suffices to justify the
claim for $P(X_1, X_2, \ldots, X_n) = X_j^m$ with $1 \le j \le n$ and
$m \ge 0$. However, this follows directly from Lemma \ref{lemma
  symmetric}. \fin

In the lack of optimal time estimates for $\F_n$ through the
probabilistic approach used so far ---see \cite{JPPP} for related
results--- we conclude this paper with optimal hypercontractivity
bounds for linear combinations of words with length lower than or
equal to 1. We will use two crucial results, the second one is
folklore and it follows from the \lq\lq invariance by rotation" of the
CAR algebra generators.
\begin{itemize}
\item \emph{The Ball}/\emph{Carlen}/\emph{Lieb convexity inequality}
  \cite{BCL}
\begin{equation}
\nonumber \Big(\frac{\mathrm{Tr} |A+B|^p + \mathrm{Tr}
    |A-B|^p}{2}\Big)^\frac{2}{p}\ge
  \big(\mathrm{Tr}|A|^p\big)^\frac{2}{p}+(p-1)\big(\mathrm{Tr}|B|^p\big)^\frac{2}{p}
\end{equation}
  for any $1 \le p \le 2$ and any given pair of $m \times m$ matrices
  $A$ and $B$.

\vskip5pt

\item \emph{A Khintchine inequality for fermion algebras} $$\Big\|
  \sum_{j=1}^d \rho_j x_j \Big\|_p \, = \, \Big( \sum_{j=1}^d
  |\rho_j|^2 \Big)^{\frac12}$$ whenever $1 \le p < \infty$,
  $\rho_j \in \R$, $x_j=x_j^*$ and $x_ix_i + x_jx_i = 2
  \delta_{ij}$.
\end{itemize}

\begin{theorem}\label{ApplicationII}
Let us denote by $\mathcal W_1$ the linear span of all words in
$\mathcal L (\F_n)$ of length lower than or equal to $1$. Then, the
following optimal hypercontractivity bounds hold for $1 < p \leq 2$,
every $t \ge -\frac12 \log (p-1)$ and all $f \in \mathcal
W_1$ $$\|\mathcal{P}_{\F_n,t} f \|_{L_2(\mathcal{L}(\F_n))} \, \le \,
\|f\|_{L_p(\mathcal{L}(\F_n))}.$$
\end{theorem}

\dem The optimality of our estimate follows once again from Weissler
argument in \cite[pp 220]{Weissler}. Moreover, it suffices to show the
inequality for the extreme case $e^{-t}=\sqrt{p-1}$. The key point in
the argument is the use of the $*$-homomorphism $\Phi:
\mathcal{L}(\F_n) \to \mathcal{L}(\mathbb{G}_{2n})$ defined in the
proof of Theorem A in conjunction with our characterization of
$\mathcal L(\mathbb{G}_{2n})$ using a spin matrix model. Indeed, we
will consider here exactly the same matrix model with $2n$ free copies
and just one generator per algebra. More precisely, given $m \ge 1$ we
will consider $x^\alpha(k)$ with $1 \le \alpha \le 2n$ and $1 \le k
\le m$ verifying the same relations as in \eqref{sign choice}
depending on the corresponding random functions $\varepsilon
((\alpha,k),(\beta,\ell))$. We also
set $$\tilde{x}^\alpha(m)=\frac{1}{\sqrt{m}}\sum_{k=1}^mx^\alpha(k)$$
as usual. Note that this model describes ---in the sense of Theorem
\ref{moments}--- the algebra $\mathcal
L(\mathbb{G}_{2n})$. In fact, according to Lemma \ref{approx norm} we
know that for every trigonometric polynomial $z = \sum_j \rho_j x_j
\in \mathcal L(\mathbb{G}_{2n})$ in the span of finite words, we can
define the corresponding elements $\tilde{z}(m) = \sum_j \rho_j
\tilde{x}_j(m) \in \mathcal A_{\varepsilon_m}$ such
that $$\lim_{m\rightarrow \infty}\|\tilde{z}(m)\|_{L_p(\mathcal
  A_{\varepsilon_m})} \, = \, \|z\|_{L_p(\mathcal
  L(\mathbb{G}_{2n}))}$$ almost everywhere. Furthermore, by dominated
convergence we find $$\hskip-12pt \lim_{m\rightarrow
  \infty}\mathbb{E}_\omega \|\tilde{z}(m)\|_{L_p(\mathcal
  A_{\varepsilon_m})} \, = \, \|z\|_{L_p(\mathcal
  L(\mathbb{G}_{2n}))}.$$ We first consider a function $f = a_0
\mathbf{1} + a_1 \lambda(g_1) + b_1 \lambda(g_1)^* + \ldots + a_n
\lambda(g_n) + b_n \lambda(g_n)^*$ in $\mathcal W_1$ such that
$\arg(a_\alpha)=\arg(b_\alpha)$ for all $1\le \alpha \le n$. By the
comments above, we have for every $1< p < 2$
\begin{eqnarray*}
\|f\|_{L_p(\mathcal L(\F_n))}^2 & = & \|\Phi f\|^2_{L_p(\mathcal
  L(\mathbb{G}_{2n}))} \\ \nonumber & = & \lim_{m \rightarrow
  \infty}\mathbb{E}_\omega \Big\| a_0 \mathbf{1} + a_1 \tilde{x}^1(m)
\tilde{x}^2(m) + b_1 \tilde{x}^2(m) \tilde{x}^1(m) \\ \nonumber & + &
\cdots \ + a_n \tilde{x}^{2n-1}(m) \tilde{x}^{2n}(m) + b_n
\tilde{x}^{2n}(m) \tilde{x}^{2n-1}(m)
\Big\|_{L_p(\mathcal{A}_{\eps_m})}^2.
\end{eqnarray*}
Now, we claim that $\|f\|_{L_p(\mathcal L(\F_n))}^2$ is bounded below
by $$\lim_{m \rightarrow \infty} \mathbb{E}_\omega \Big(|a_0|^2 +
\frac{p-1}{m^2} \sum_{\begin{subarray}{c} 1\le \alpha \le n \\ 1\le k
    \le m \end{subarray}} \Big\|\sum_{1 \le \ell \le m} \Big( a_\alpha
+ b_\alpha \varepsilon
\big((2\alpha-1,k),(2\alpha,\ell)\big)\Big)x^{2\alpha}(\ell)
\Big\|_p^2 \Big).$$ If this is true, we can apply Khintchine's
inequality for fixed $\alpha$ and $k$ to get
\begin{eqnarray*}
\lefteqn{\mathbb{E}_\omega \Big\|\sum_{1\leq \ell \leq m}
  \Big( a_\alpha + b_\alpha
  \varepsilon\big((2\alpha-1,k),(2\alpha,\ell) \big)\Big)
  x^{2\alpha}(\ell) \Big\|_p^2} \\ & = & \mathbb{E}_\omega
\Big\|\sum_{1\leq \ell \leq m} \Big( |a_\alpha| + |b_\alpha|
\varepsilon\big((2\alpha-1,k),(2\alpha,\ell) \big)\Big)
x^{2\alpha}(\ell) \Big\|_p^2 \\ & =
&\sum_{1\le \ell \le m} \Big( |a_\alpha|^2 + |b_\alpha|^2 + 2
|a_\alpha b_\alpha| \mathbb{E}_\omega
\varepsilon\big((2\alpha-1,k),(2\alpha,\ell) \big) \Big) \, =
\, m(|a_\alpha|^2+ |b_\alpha|^2).
\end{eqnarray*}
Here, we have used that the $\varepsilon$'s are centered for $\alpha
\neq \beta$. Therefore, we finally obtain $$\|f\|_{L_p(\mathcal
  L(\F_n))}^2 \, \ge \, |a_0|^2 + (p-1)
\sum_{\alpha=1}^n(|a_\alpha|^2+|b_\alpha|^2) \, = \,
\|\mathcal{P}_{\F_n,t}f\|_{L_2(\mathcal{L}(\F_n))}^2$$ for
$e^{-t}=\sqrt{p-1}$. Therefore, it suffices to prove the claim. To
this end, note that $$\|f\|_{L_p(\mathcal L(\F_n))}^2 \, = \, \lim_{m
  \rightarrow \infty} \mathbb{E}_\omega \big\|A_m + x^1(1) B_m
\big\|_{L_p(\mathcal{A}_{\eps_m})}^2,$$ where $A_m$ and $B_m$ are
given by
\begin{eqnarray*}
A_m & = & a_0 \mathbf{1} + \frac{1}{m} \sum_{\substack{2\le k \le m
    \\ 1 \le \ell \le m}} \Big( a_1 + b_1 \varepsilon\big( (1,k),
(2,\ell) \big) \Big) x^1(k)x^2(\ell) \\ & + & \frac{1}{m} \sum_{1 \le
  k,\ell \le m} \Big[ a_2 x^3(k) x^4(\ell) + b_2 x^4(k) x^3(\ell) +
  \ldots + b_n x^{2n}(k) x^{2n-1}(\ell) \Big]
\end{eqnarray*}
and $B_m = \frac{1}{m} \sum_{1 \leq \ell \leq m} \big( a_1 + b_1
\varepsilon\big((1,1),(2,\ell) \big) \big) x^2(\ell)$. Then, since the
spin matrix model is unaffected by the change of sign of one generator
and $A_m,B_m$ do not depend on $x^1(1)$, we deduce
$\|A_m+x^1(1)B_m\|_p=\|A_m-x^1(1)B_m\|_p$. Therefore, applying
Ball/Carlen/Lieb inequality we conclude that $$\|f\|_{L_p(\mathcal
  L(\F_n))}^2 \, \ge \, \lim_{m\rightarrow \infty}\mathbb{E}_\omega
\Big( \|A_m\|_{L_p(\mathcal{A}_{\eps_m})}^2 + (p-1)
\|B_m\|_{L_p(\mathcal{A}_{\eps_m})}^2 \Big),$$ where we have used that
$\|x^1(1)B_m\|_p=\|B_m\|_p$ for every $\omega$ and every $p$. If we
apply the same strategy with $x^1(2),\ldots, x^1(m)$, it is not
difficult to obtain the following lower bound
\begin{eqnarray*}
\|f\|_{L_p(\mathcal L(\F_d))}^2 & \ge & \lim_{m \rightarrow \infty}
\mathbb{E}_\omega \Big\| a_0 \mathbf{1} + a_2 \tilde{x}^3(m)
\tilde{x}^4(m) + b_2 \tilde{x}^4(m) \tilde{x}^3(m) \\ [3pt] & + &
\cdots \ + a_n \tilde{x}^{2n-1}(m) \tilde{x}^{2n}(m) + b_n
\tilde{x}^{2n}(m) \tilde{x}^{2n-1}(m) \Big\|_p^2 \\ & + &
\frac{p-1}{m^2} \sum_{1 \le k \le m} \Big\| \sum_{1 \le \ell \le
  m}\Big( a_1 + b_1 \varepsilon\big((1,k),(2,\ell)\big) \Big)x^2(\ell)
\Big\|_p^2.
\end{eqnarray*}
Our claim follows iterating this argument on $2 \le \alpha \le n$.  It
remains to consider an arbitrary $f = a_0 \mathbf{1} + a_1
\lambda(g_1) + b_1 \lambda(g_1)^* + \ldots + a_n \lambda(g_n) + b_n
\lambda(g_n)^* \in \mathcal W_1$.  Let us set
$(\theta_\alpha, \theta'_\alpha) = (\arg(a_\alpha), \arg(b_\alpha))$
and $(\nu_\alpha, \nu'_\alpha) = (\frac12 (\theta_\alpha+\theta'_\alpha),
\frac12 (\theta_\alpha-\theta'_\alpha))$ for each $1\le \alpha \le n$. Consider the
1-dimensional representation $\pi: \F_n \to \C$ determined by
$\pi(g_\alpha) = \exp(i\nu'_\alpha )$ for the $\alpha$-th generator
$g_\alpha$. According to the $L_p$-analog of Fell's absorption
principle \cite{PP}, we have from the first part of the proof that
\begin{eqnarray*}
\|\mathcal{P}_{\F_n,t} f\|_2 & \le & \Big\| a_0 \mathbf{1} +
\sum_{\alpha=1}^n |a_\alpha|e^{i\nu_\alpha} \lambda(g_\alpha) +
|b_\alpha|e^{i\nu_\alpha} \lambda(g_\alpha)^*
\Big\|_{L_p(\mathcal{L}(\F_n))} \\ & = & \Big\| a_0 \mathbf{1} +
\sum_{\alpha=1}^n |a_\alpha|e^{i\nu_\alpha} \pi(g_\alpha)
\lambda(g_\alpha) + |b_\alpha|e^{i\nu_\alpha}
\pi(g_\alpha^{-1})\lambda(g_\alpha)^* \Big\|_{L_p(\mathcal{L}(\F_n))}
\\ & = & \Big\| a_0 \mathbf{1} + \sum_{\alpha=1}^n a_\alpha
\lambda(g_\alpha) + b_\alpha \lambda(g_\alpha)^*
\Big\|_{L_p(\mathcal{L}(\F_n))} \ = \ \|f\|_{L_p(\mathcal{L}(\F_n))}.
\end{eqnarray*}
The proof is complete. \fin

We finish this section with further results on $L_p\to L_2$ estimates for the free Poisson semigroup.
The key point here is to use a different model for Haar unitaries.
In the sequel, we will denote by $\mathbb M_2$ the algebra of $2\times 2$ matrices.
\begin{lemma}\label{haar}
If $u_j=\lambda(g_j)$ and $x_j=\lambda(z_j)$, the map $$u_j \mapsto \left[\begin{array}{cc} 0 & x_{2j-1} \\ x_{2j}&
0\end{array}\right]$$ determines a trace preserving $*$-homomorphism
$\pi: \mathcal L(\mathbb F_n) \to \mathbb M_2\otens \mathcal L(\mathbb
G_{2n})$ such that $$\pi \circ \mathcal P_{\F_n, t} \, = \, \big( Id_{\mathbb M_2} \otimes \mathcal P_{\mathbb{G}_{2n}, t} \big) \circ \pi.$$
\end{lemma}

\dem Since $\pi(u_j)$ is a unitary $w_j$ in $\mathbb M_2\otens \mathcal L(\mathbb{G}_{2n})$ and $\mathbb F_n$ is a free group, a unique $*$-homomorphism $\pi: \mathcal L(\mathbb F_n) \to \mathbb M_2\otens \mathcal L(\mathbb{G}_{2n})$ is determined by the $w_j$'s. Thus, it suffices to check that $\pi$ is trace preserving. The fact that $\pi(\lambda(g))$ has trace zero in $\mathbb M_2\otens \mathcal L(\mathbb G_{2n})$ for every $g \neq e$ follows easily from the equalities $$\pi(u_1)^{2k}= \left[\begin{array}{cc} a_{1,2k} &0 \\ 0 & b_{1,2k}\end{array}\right] \qquad \pi(u_1)^{2k+1} = \left[\begin{array}{cc} 0& a_{1,2k+1} \\ b_{1,2k+1} &
0 \end{array}\right]$$ and its analogous formulae for the product of different generators. Here, we have used the notations introduced in \eqref{defab}. The second assertion can be checked by simple calculations. The proof is complete. \fin
Biane's theorem relies on an induction argument \cite[Lemma 2]{Biane} which exploits the Ball-Carlen-Lieb convexity inequality stated before Theorem \ref{ApplicationII}. In fact, our proof of Theorem \ref{ApplicationII} follows the same induction argument. We will now consider spin matrix models with operator coefficients. More precisely, given a finite von Neumann algebra $(\mathcal M,\tau)$, we will look at $\mathcal M\otens \mathcal A_{\varepsilon_m}$. In particular, following the notation in Section \ref{Section: proof of theorem B} every $x\in\mathcal M\otens \mathcal A_{\varepsilon_m}$ can be written as $x=\sum_{A} \rho_A \tens x_A^{\varepsilon_m}$ where $\rho_A \in \mathcal M$ for every $A \subset\Upsilon_m$. Then, the induction argument easily leads to the inequality below provided that $e^{-t}\leq \sqrt{p-1}$
$$\|x\|^2_{L_p(\mathcal M\otens \A_{\varepsilon_m})} \geq
\sum_{A\subset\Upsilon_m} e^{-2t|A|} \| \rho_A \|_{L_p(\mathcal M)}^2.$$
For our purpose we will consider $\mathcal M=\mathbb M_2$ with its normalized trace, so that $$\|a\|_p \geq 2^{\frac 1 2 - \frac 1 p} \|a\|_2$$ for every $a \in \mathbb  M_2$. Let $x=\sum_{A} \rho_A \tens x_A^{\varepsilon_m}$ be as above.  Let us also define $\mathcal U$ as the (possible empty) set of the subsets $A$ of $\Upsilon_m$ such that $\rho_A$ is a multiple of a unitary. In particular, $\|\rho_A\|_{L_2(\mathcal M)}=\|\rho_A\|_{L_p(\mathcal M)}$ for every $A\in \mathcal U$. Then, letting $y=\sum_{A\in \mathcal U} \rho_A \tens x_A^{\varepsilon_m}$, the following estimate holds provided $e^{-t}\leq \sqrt{p-1}$
\begin{equation}\label{Bianecoef} \| x \|^2_{p} \geq \big\| Id_{\mathbb
  M_2}\tens \mathcal{S}_{\varepsilon_m,t}(y)\big\|^2_{2} +
  2^{1-\frac 2 p } \big\| Id_{\mathbb M_2}\tens
  \mathcal{S}_{\varepsilon_m,t}(x-y)\big\|^2_{2},
\end{equation}
where the right hand side norms are taken in $\mathbb{M}_2\otens \A_{\varepsilon_m}$. Our first application of this alternative approach is that Weissler's theorem \cite{Weissler} can be proved using probability and operator algebra methods.
\begin{proposition}\label{Weissler}
If $1 < p \le q < \infty$, we find $$\big\| \mathcal{P}_{\mathbb{Z},t}: L_p(\mathcal{L}(\mathbb{Z})) \to
L_q(\mathcal{L}(\mathbb{Z})) \big\| \, = \, 1 \ \Leftrightarrow \ t \ge \frac12 \log \frac{q-1}{p-1}.$$
\end{proposition}

\dem We will assume that $q=2$ since the optimal time for every $p,q$ can be obtained from this case by means of standard arguments involving log-Sobolev inequalities. We follow here the same approximation procedure of Lemmas \ref{approx norm} and \ref{key lemma} with $n=2$ and $d=1$. Consider a reduced word $x=x_{\alpha_1}\cdots x_{\alpha_s}$ in
$\mathcal L(\mathbb G_2)$, so that $\alpha_j \in\{1,2\}$ and $\alpha_j\neq \alpha_{j+1}$. We then form the associated element
$$\tilde{x}(m)(\omega) = \frac{1}{m^{s/2}} \sum_{\substack{\k \in
    [m]^s \\ \sigma(\k) = \sigma_0}} x^{\alpha_1}(k_1)(\omega) \cdots
x^{\alpha_s}(k_s)(\omega)\in \A_{\varepsilon_m}
.$$ Note that restricting to $\sigma(\k) = \sigma_0$ implies that there will be no repetitions of the elements $x^{\alpha_j}(k_j)$, hence no simplifications in $\tilde{x}(m)$. As we showed in the proof of Lemma \ref{key lemma}, the terms with repetitions do not play any role. On the other hand, Lemma \ref{approx norm} easily extends to operator coefficients so that for any $1\leq p\leq 2$, every $\rho_j \in \mathbb{M}_2$ and every reduced word $x_j\in \mathcal L(\mathbb G_2)$, we have
\begin{equation}\label{normcoef}
\lim_{m\to \infty} \Big\| \sum_j \rho_j \tens \tilde x_j(m)\Big\|_{L_p(\mathbb M_2\otens \A_{\varepsilon_m})}
=\Big\|\sum_j \rho_j \tens x_j\Big\|_{L_p(\mathbb M_2\otens \mathcal L(\mathbb G_2))} \quad
\mbox{a.e.}
\end{equation}
Let us denote by $u=\lambda(g_1)$ the canonical generator of $\mathcal L (\mathbb Z)$.  By the positivity of $\mathcal P_{\Z,t}$ and a density argument, it suffices to show that $\|\mathcal P_{\Z,t} f\|_{L_2(\mathcal L (\mathbb Z))}\leq \|f\|_{L_p(\mathcal L (\mathbb Z))}$ for every positive trigonometric polynomial $$\displaystyle{f=\rho_0 \mathbf{1}+\sum_{j=1}^d
(\rho_j u^j+\overline \rho_j u^{*j})}.$$ To this end, we use the map $\pi$ from Lemma \ref{haar} and construct
\begin{eqnarray*}
x \ = \ \pi (f) & = & \left[\begin{array}{cc} \rho_0 & 0 \\ 0&
    \rho_0\end{array}\right] \tens \mathbf{1} \\ & + & \displaystyle{\sum_{\ell \ge 1}}
\left[\begin{array}{cc} \rho_{2\ell} & 0 \\ 0& \overline
    \rho_{2\ell}\end{array}\right]\tens a_{1,2\ell} +
\left[\begin{array}{cc} \overline \rho_{2\ell} & 0 \\ 0&
    \rho_{2\ell}\end{array}\right]\tens b_{1,2\ell} \\ & + &
\displaystyle{\sum_{\ell \ge 1}} \left[\begin{array}{cc} 0 &
    \rho_{2\ell+1} \\ \overline \rho_{2\ell+1} &
    0\end{array}\right]\tens a_{1,2\ell+1} +
\left[\begin{array}{cc}0 & \overline \rho_{2\ell+1} \\ \rho_{2\ell+1}&
    0\end{array}\right]\tens b_{1,2\ell+1}.
\end{eqnarray*}
To use our approximation procedure, we consider the element $\tilde x(m)\in \mathbb M_2\otens \A_{\varepsilon_m}$ associated to
$x$. We start noting that $\tilde x(m)$ is self-adjoint. Now, in order to use \eqref{Bianecoef} and make act $Id_{\mathbb M_2}\otimes \mathcal S_{\varepsilon_m,t}$, we must write $\tilde x(m)$ in reduced form. That is, for every $\k \in [m]^s$ with $\sigma(\k) = \sigma_0$ and $\underline{\alpha}=(\alpha_1,\cdots,\alpha_s)\in \{1,2\}^s$ with $\alpha_j\neq \alpha_{j+1}$, we want to understand the matrix coefficients $\gamma^{\underline{\alpha}}(\k)$ of  $x^{\underline{\alpha}}(\k)=x^{\alpha_1}(k_1) \cdots x^{\alpha_s}(k_s)$, where the latter is an element in the basis of $\A_{\varepsilon_m}$. In fact, it suffices to show that these matrix coefficients are multiples of unitaries, so that all the subsets $A$ of $\Upsilon_m$ are in $\mathcal U$ and we do not loose any constant when applying \eqref{Bianecoef}.
Let us first assume that $s=2\ell+1$ is odd.
Since by definition there is no simplifications in $\tilde x(m)$, the term $x^{\underline{\alpha}}(\k)$ will only appear in the element in $\A_{\varepsilon_m}$ associated to either $a_{1,2\ell+1}$ or $b_{1,2\ell+1}$.
By the commutation relations, we see that $x^{\underline{\alpha}}(\k)^*=\pm x^{\underline{\alpha}}(\k)$. Then its matrix coefficient must also satisfy $\gamma^{\underline{\alpha}}(\k)^*=\pm \gamma^{\underline{\alpha}}(\k)$.
Moreover, one easily checks that it also has the shape $$\left[\begin{array}{cc} 0 & \delta \\ \mu
    &0\end{array}\right]$$ from the above formula of $x$. Hence $\delta=\pm
\overline \mu$ and $\gamma^{\underline{\alpha}}(\k)$ is a multiple of a unitary (this
can also be directly seen from the formula of $x$).
If $s=2\ell$, the term $x^{\underline{\alpha}}(\k)$ will appear in the elements associated to the two reduced words $a_{1,2\ell}$ and $b_{1,2\ell}$.
Since the commutation relations only involve signs, after a moment of thought we
can conclude that $\gamma^{\underline{\alpha}}(\k)$ has the shape
$$\left[\begin{array}{cc} \delta& 0 \\ 0 &\overline
    \delta\end{array}\right].$$ Hence, it is a multiple of a
unitary. Actually, we also know that $\delta$ is either real or purely imaginary. Once we have seen that the matrix coefficients of $\tilde x(m)$ written in reduced form are multiples of unitaries, we can conclude the proof as in Theorem B. Indeed, using Lemma \ref{haar}, \eqref{Bianecoef} and \eqref{normcoef}, we get
\begin{eqnarray*}
\| f\|_{L_p(\mathbb T)} & = & \| x\|_{L_p(\mathbb M_2\otens \mathcal L (\mathbb
G_2))} \\ &=& \lim_{m\to \infty} \|\tilde x(m)\|_{L_p(\mathbb M_2\otens \A_{\varepsilon_m})} \\ &\geq& \lim_{m\to
   \infty} \|(Id_{\mathbb M_2}\tens \mathcal S_{\varepsilon_m,t})\tilde x(m)\|_{L_2(\mathbb M_2\otens \A_{\varepsilon_m})} \\&= & \| (Id_{\mathbb M_2} \tens
 \mathcal P_{\mathbb{G}_2,t}) (x)\|_{L_2(\mathbb M_2\otens \mathcal L(\mathbb
G_2))} =\|\mathcal P_{\Z,t}(f)\|_{L_2(\mathbb T)},\end{eqnarray*}where the limits are taken a.e. and $t\geq -\frac 1 2 \log(p-1)$. The proof is complete. \fin
A slight modification of the previous argument allows us to improve Theorem A ii) for $q=2$. In fact, by a standard use of log-Sobolev inequalities we may also improve the $L_p \to L_q$ hypercontractivity bound, see Remark \ref{pqHC} below.
\begin{theorem}\label{improvement A iii)}
If $1< p\leq 2$, we find $$\big\|\mathcal  P_{\mathbb F_n, t} : L_p(\mathcal L (\mathbb F_n))\to L_2(\mathcal L (\mathbb F_n))\big\|=1 \quad {if }\quad t\geq \frac 1 2 \log \frac{1}{p-1} + \frac 1 2 \Big(\frac 1
 p - \frac 1 2\Big) \log 2 .$$
\end{theorem}

\dem Once again, by positivity and density it suffices to prove the assertion for a
positive trigonometric polynomial $f\in \mathcal L (\mathbb F_n)$. If $\jbar=(j_1, \ldots, j_d)$, we will use the notation $|\jbar|=d$ and $u_{\jbar}=\lambda(g_{\jbar})$ with $g_{\jbar}=g_{j_1}\cdots g_{j_d}$ a reduced word in $\F_n$, so that $$f=\summ_{\jbar} \rho_{\jbar} u_{\jbar}.$$ Here we use the usual convention that $g_{-k}=g_k^{-1}$. We use again the trace preserving $*$-homomorphism $\pi: \mathcal L
(\mathbb F_n)\to \mathbb M_2 \otens \mathcal L (\mathbb G_{2n})$
coming from Lemma \ref{haar}. This gives the identity
$$\pi(u_{\jbar})=\left[\begin{array}{cc} 0 & x_{2j_1-1}
    \\ x_{2j_1}& 0\end{array}\right] \left[\begin{array}{cc} 0 & x_{2j_2-1}
    \\ x_{2j_2}& 0\end{array}\right] \cdots\left[\begin{array}{cc} 0 & x_{2j_d-1}
    \\ x_{2j_d}& 0\end{array}\right]$$
with the convention that for $j>0$, $x_{-2j}=x_{2j-1}$ and $x_{-2j-1}=x_{2j}$. If $d=0$, we set $g_{\jbar} =e$ and $\pi(u_{\jbar})=Id_{\mathbb M_2}$.
Hence with $x=\pi(f)$, summing up according to the length we obtain
\begin{equation} \label{formula x=pi(f)}
\begin{array}{rclcl} x & \!\!\! = \!\!\! & \left[\begin{array}{cc} \rho_0 & 0 \\ 0& \rho_0 \end{array} \right]
\tens \mathbf{1} & \!\!\! + \!\!\! & \displaystyle{\sum_{\substack{|\jbar|=2\ell \\ \ell \geq 1}}}
 \left[\begin{array}{cc} \rho_{\jbar} & 0 \\ 0&
    \rho_{-\jbar}\end{array}\right] \tens x_{2j_1-1}x_{2j_2} \cdots x_{2j_{2\ell}}  \\
& & & \!\!\! + \!\!\! & \displaystyle{\sum_{\substack{|\jbar|=2\ell+1 \\ \ell \geq 0}}}
 \left[\begin{array}{cc} 0 &\rho_{\jbar} \\
    \rho_{-\jbar}&0\end{array}\right] \tens x_{2j_1-1}x_{2j_2} \cdots x_{2j_{2\ell+1}-1}. \end{array}
\end{equation}
We repeat the arguments used in the proof of Proposition \ref{Weissler} to approximate
$x$ by a spin model $\tilde x(m)(\omega)$ with operator coefficients. That is, $x_{\alpha_1}\cdots x_{\alpha_s}\in \mathcal{L}(\mathbb G_{2n})$ is associated to $$\tilde{x}(m)(\omega) = \frac{1}{m^{s/2}} \sum_{\substack{\k \in
    [m]^s \\ \sigma(\k) = \sigma_0}} x^{\alpha_1}(k_1)(\omega) \cdots
x^{\alpha_s}(k_s)(\omega)\in \A_{\varepsilon_m}.$$
Note that the contribution to $x$ given by \eqref{formula x=pi(f)} of words of length 0 and 1 is
\begin{eqnarray*}
\left[\begin{array}{cc} \rho_0 & 0 \\ 0& \rho_0\end{array}\right]
\tens 1 +\displaystyle{\sum_{j\in \mathbb Z \setminus \{0\}}}
\left[\begin{array}{cc} 0 & \rho_{j} \\ \rho_{-j} &
    0\end{array}\right]\tens x_{2j-1}.  \end{eqnarray*}
Since $f$ is self-adjoint, we have $\rho_{-j}=\overline \rho_j$ for $j\in \Z \setminus \{0\}$.
Hence the matrix coefficients corresponding to the words of length 0 and 1 in the approximation are multiples of unitaries.
We will have $\{A\subset \Upsilon_m \; :\; |A|\leq 1\}\subset \mathcal U$ with the notations of \eqref{Bianecoef},
and decompose $f=g+h$, where $g$ is the part of $f$ of degree less than 1 and $h$ is supported by the words of length greater or equal than 2.
Observe that $g$ and $h$ are orthogonal.
Let $t=t_0+t_1$ with $t_0=-\frac 1 2 \log(p-1)$.
Since $h$ has valuation 2, we have
$$\|\mathcal P_{\F_n,t_0+t_1}(h)\|_{2}\leq e^{-2t_1} \|\mathcal P_{\F_n,t_0}(h)\|_{2}.$$
Thus thanks to
\eqref{Bianecoef}, as in the proof of Proposition \ref{Weissler}, we get by orthogonality
\begin{eqnarray*}
\|f\|_{p}^2 & \geq & \| \mathcal P_{\F_n,t_0} (g)\|_{2}^2 + 2^{1-\frac 2p}
\|\mathcal P_{\F_n,t_0} (h)\|_{2}^2 \\ & \geq & \| \mathcal P_{\F_n,t} (g)\|_{2}^2 + 2^{1-\frac 2p}e^{4t_1}
\|\mathcal P_{\F_n,t} (h)\|_{2}^2 \\ & \geq & \| \mathcal P_{\F_n,t} (g)\|_{2}^2 + \|\mathcal P_{\F_n,t} (h)\|_{2}^2
=  \| \mathcal P_{\F_n,t} (f)\|_{2}^2,
\end{eqnarray*}
provided that $e^{-4t_1} 2 ^{\frac 2 p-1} \leq 1 \Leftrightarrow t_1 \ge \frac 1 2 (\frac 1
 p - \frac 1 2) \log 2$. This completes the proof. \fin

\begin{remark} \rm{
Let $\sigma$ be the involutive $*$-representation on $\mathcal L(\mathbb F_n)$
exchanging $u_j$ and $u_j^*=u_{-j}$ for all $j\geq 1$. So that if
$f=\sum_{\jbar} \rho_{\jbar} u_{\jbar}$, then
$\sigma(f)=\sum_{\jbar} \rho_{-\jbar} u_{\jbar}$. Denote by $\mathcal
L(\mathbb F_n)^\sigma$ the fixed point algebra of $\sigma$, it clearly contains
$\mathcal A_{sym}^n$. The above arguments actually prove
that $\mathcal P_{\F_n,t}$ is hypercontractive on $\mathcal
L(\mathbb F_n)^\sigma$ from $L_p$ to $L_2$ with optimal time.
Indeed, under this symmetric condition for $f$ all the matrix coefficients will be multiples of unitaries.
Then using
the equivalence between hypercontractivity with optimal time and log-Sobolev inequality, one sees that Theorem A iii) can be extended to $\mathcal
L(\mathbb F_n)^\sigma$. The Gross' argument to deduce general hypercontractive inequalities $L_p\to L_q$ from $L_p\to L_2$ estimates in this setting are recalled in \cite{JPPP}.}
\end{remark}

\begin{remark} \label{pqHC} \rm{
We claim that
$$ \frac 1 2 \log \frac{1}{p-1} + \frac 1 2 \Big(\frac 1 p - \frac 1 2\Big) \log 2
\leq \frac{\beta}{2} \log \frac{1}{p-1}$$
with $\beta=1+\frac{\log(2)}{4}$.
This is not difficult to prove by using basic computations.
Then in particular Theorem \ref{improvement A iii)} proves that we have hypercontractive $L_p\to L_2$ estimates for
$t\geq -\frac{\beta}{2} \log(p-1)$.
Then Gross' arguments relying on log-Sobolev inequality apply when the time has this shape, and the constant $2$ given by Theorem A ii) can be replaced by the better constant $\beta=1+\frac{1}{4} \log(2) \sim 1.17$. Hence for any $1 < p \le q < \infty$ we get
$$\big\| \mathcal{P}_{\F_n,t} \hskip1pt : L_p( \hskip0.5pt
  \mathcal{L}(\F_n) \hskip0.5pt ) \to L_q( \hskip0.5pt
  \mathcal{L}(\F_n) \hskip0.5pt ) \hskip1pt \big\| \, = \, 1
  \hskip13pt \mbox{if} \hskip13pt t \ge \frac{\beta}{2} \log \frac{q-1}{p-1}.$$
}
\end{remark}

\bibliographystyle{amsplain}


\

\hfill {\bf Marius Junge} \\ \null \hfill Department of Mathematics
\\ \null \hfill University of Illinois at Urbana-Champaign \\ \null
\hfill 1409 W. Green St. Urbana, IL 61891. USA \\ \null \hfill
\texttt{junge@math.uiuc.edu}

\vskip5pt


\hfill \noindent \textbf{Carlos Palazuelos} \\ \null \hfill Instituto
de Ciencias Matem{\'a}ticas \\ \null \hfill CSIC-UAM-UC3M-UCM \\ \null
\hfill Consejo Superior de Investigaciones Cient\'{\i}ficas \\ \null
\hfill C/ Nicol\'as Cabrera 13-15.  28049, Madrid. Spain \\ \null
\hfill\texttt{carlospalazuelos@icmat.es}

\vskip5pt

\hfill \noindent \textbf{Javier Parcet} \\ \null \hfill Instituto de
Ciencias Matem{\'a}ticas \\ \null \hfill CSIC-UAM-UC3M-UCM \\ \null
\hfill Consejo Superior de Investigaciones Cient\'{\i}ficas \\ \null
\hfill C/ Nicol\'as Cabrera 13-15.  28049, Madrid. Spain \\ \null
\hfill\texttt{javier.parcet@icmat.es}

\vskip5pt

\hfill \noindent \textbf{Mathilde Perrin} \\ \null \hfill Instituto de
Ciencias Matem{\'a}ticas \\ \null \hfill CSIC-UAM-UC3M-UCM \\ \null
\hfill Consejo Superior de Investigaciones Cient\'{\i}ficas \\ \null
\hfill C/ Nicol\'as Cabrera 13-15.  28049, Madrid. Spain \\ \null
\hfill\texttt{mathilde.perrin@icmat.es}

\vskip5pt

\hfill \noindent \textbf{\'Eric Ricard} \\ \null \hfill Laboratoire de
Math\'ematiques Nicolas Oresme \\ \null \hfill Universit\'e de Caen
Basse-Normandie \\ \null \hfill 14032 Caen Cedex, France \\ \null
\hfill\texttt{eric.ricard@unicaen.fr}

\end{document}